\documentclass[11pt,a4paper]{article}
\usepackage[margin=1in]{geometry}
\usepackage{latexsym, amsfonts, amsmath}
\newcommand{\be}{\begin{equation}}
\newcommand{\ee}{\end{equation}}
\newcommand{\bea}{\begin{eqnarray}}
\newcommand{\eea}{\end{eqnarray}}
\newcommand{\bean}{\begin{eqnarray*}} 
\newcommand{\eean}{\end{eqnarray*}}

\newcommand{\brray}{\begin{array}}
\newcommand{\erray}{\end{array}}
\newcommand{\ben}{\begin{equation}{nonumber}}
\newcommand{\een}{\end{equation}{nonumber}}

\newtheorem{dfn}{Definition}[section]
\newtheorem{thm}[dfn]{Theorem}
\newtheorem{lmma}[dfn]{Lemma}
\newtheorem{ppsn}[dfn]{Proposition}
\newtheorem{crlre}[dfn]{Corollary}
\newtheorem{xmpl}[dfn]{Example}
\newtheorem{rmrk}[dfn]{Remark}
\newcommand{\bdfn}{\begin{dfn}}
\newcommand{\bthm}{\begin{thm}}
\newcommand{\blmma}{\begin{lmma}}
\newcommand{\bppsn}{\begin{ppsn}}
\newcommand{\bcrlre}{\begin{crlre}}
\newcommand{\bxmpl}{\begin{xmpl}}
\newcommand{\brmrk}{\begin{rmrk}}
\newcommand{\edfn}{\end{dfn}}
\newcommand{\ethm}{\end{thm}}
\newcommand{\elmma}{\end{lmma}}
\newcommand{\eppsn}{\end{ppsn}}
\newcommand{\ecrlre}{\end{crlre}}
\newcommand{\exmpl}{\end{xmpl}}
\newcommand{\ermrk}{\end{rmrk}}

\newcommand{\IC}{\mathbb{C}}

\newcommand{\IR}{\mathbb{R}}

\newcommand{\cla}{{\cal A}}
\newcommand{\clb}{{\cal B}}
\newcommand{\clc}{{\cal C}}
\newcommand{\cld}{{\cal D}}
\newcommand{\cle}{{\cal E}}
\newcommand{\clf}{{\cal F}}

\newcommand{\clh}{{\cal H}}
\newcommand{\cli}{{\cal I}}

\newcommand{\clk}{{\cal K}}
\newcommand{\cll}{{\cal L}}
\newcommand{\clm}{{\cal M}}
\newcommand{\cln}{{\cal N}}
\newcommand{\clo}{{\cal O}}
\newcommand{\clp}{{\cal P}}
\newcommand{\clq}{{\cal Q}}

\newcommand{\cls}{{\cal S}}

\newcommand{\clu}{{\cal U}}

\newcommand{\clw}{{\cal W}}

\def\a*{{\cal A}_{h,*}}
\def\B{{\cal B}(h)}
\def\B1{{\cal B}_1(h)}
\def\b{{\cal B}^{\rm s.a.}(h)}
\def\b1{{\cal B}^{\rm s.a.}_1(h)}

\newcommand{\ot}{\otimes}

\newcommand{\raro}{\rightarrow}

\newcommand{\id}{\mbox{id}}

\def \qed {$\Box$}

\def\a*{{\cal A}_{h,*}}
\def\B{{\cal B}(h)}
\def\B1{{\cal B}_1(h)}
\def\b{{\cal B}^{\rm s.a.}(h)}
\def\b1{{\cal B}^{\rm s.a.}_1(h)}

\begin{document}
\begin{center}
{\Large{\bf Smooth actions of compact quantum groups on compact smooth manifolds}}\\ 
{\large {\bf Debashish Goswami\footnote{Partially supported by Swarnajayanti Fellowship from D.S.T. (Govt. of India) and also acknowledges
 the Fields Institute, Toronto for providing hospitality for a brief stay when a small part of this work was done.} and \bf Soumalya Joardar \footnote{Acknowledges support from CSIR.}}}\\ 
Indian Statistical Institute\\
203, B. T. Road, Kolkata 700108\\
Email: goswamid@isical.ac.in,\\
Phone: 0091 33 25753420, Fax: 0091 33 25773071\\
\end{center}
\begin{abstract}
Definition of a smooth action of a CQG on a compact, smooth manifold is given and studied. It is shown that a smooth action is always injective. Furhtermore A necessary and sufficient condition for a lift of the smooth action as a bimodule morphism on the bimodule of one forms has been deduced and it is also shown to be equivalent to the condition of preserving some Riemannian inner product on the manifold.  

\end{abstract}
{\bf Subject classification :} 81R50, 81R60, 20G42, 58B34.\\  
{\bf Keywords:} Compact quantum group, Riemannian manifold, smooth action.
 \section{Introduction}
 It is  a very important and interesting problem in the theory of quantum groups and noncommutative geometry to study `quantum symmetries' of various classical and quantum structures. 
 Indeed, symmetries of physical systems (classical or quantum) were conventionally modeled by group actions, and after the advent of quantum groups, group symmetries were naturally generalized to  symmetries given by quantum group action. In the framework of Connes' non commutative geometry it is  natural to consider actions of compact quantum groups on spectral triples. In the topological setting the action of a compact quantum group (CQG for short) on a spectral triple was defined in \cite{Goswami} as the $C^{\ast}$ action of the compact quantum group on the natural $C^{\ast}$ algebra in the spectral triple. For the classical case it is nothing but a continuous action of the CQG on $C(M)$ where $M$ is a compact manifold. It is worth mentioning that $C^*$ action of a CQG on a $C^*$ algebra has been extensively studied in \cite{Podles}, \cite{Act_com} and in other places. But in the context of noncommutative geometry the `space' has some additional structures and thus one expects to go beyond the $C^*$ 
action category. As in the group case, we should be able to talk about `smooth' action of a CQG. In this paper one of our jobs has been to define and study a `smooth' action of a CQG on a classical spectral triple, i.e. we consider the topological action of a CQG on the smooth algebra $C^{\infty}(M)$, where the algebra is endowed with its canonical Fr\'echet topology coming from derivations. We proved an interesting result about injectivity of a smooth action.\\
\indent For a smooth action of a group on a compact smooth manifold, the `differential' of the action automatically lifts as a well defined bimodule morphism to the space of one forms of the manifold. But in case of CQG it turns out that this lift is not automatic due to a fundamental noncommutativity. In fact, we have example of Hopf-algebra (of non compact type) having coaction on a coordinate algebra of an algberaic variety which does not admit such a lift. We give a necessary and sufficient condition for such a lift to exist. However, no such example is there yet with a smooth action of a compact quantum group.\\
\indent We show that a smooth action is inner product preserving with respect to some Riemannian metric on the manifold if and only if it admits a lift $\beta:\Omega^{1}(C^{\infty}(M))\raro \Omega^{1}(C^{\infty}(M))\bar{\ot}\clq$ to the bimodule of one forms. Already the lift played a crucial role in studying the isometric sction of a CQG on a compact, connected, Riemannian manifold. We believe that in the context of noncommutative geometry (\cite{con}) such lifts are going to be important when we are going to treat CQG as `symmetry' objects. More specifically, when we are going to consider invariance of various action functionals (Yang-Mills, Einstein-Hilbert) under CQG actions, such lifts might be very important.

\section{ Preliminaries}
In this paper all the
Hilbert spaces are over $\mathbb{C}$ unless mentioned otherwise. If $V$ is a
vector space over real numbers we denote its complexification by
$V_{\mathbb{C}}$. For a vector space $V$, $V^{'}$ stands for its algebraic dual.
$\oplus$ and $\ot$ will denote the algebraic direct sum and
algebraic tensor
product respectively. We shall denote the $C^{\ast}$
algebra of bounded operators on a Hilbert space $\clh$ by $\clb(\clh)$ and the
$C^{\ast}$ algebra
of compact operators on $\clh$ by $\clb_{0}(\clh)$. $Sp$, $\overline{Sp}$)
stand for the
linear span and closed linear span of elements of a vector space respectively, whereas ${\rm Im}(A)$ denotes the image of a linear map. 
We denote by  WOT and SOT  the weak operator
topology and the strong operator topology respectively. Let $\clc$ be an
algebra. Then
$\sigma_{ij}:\underbrace{\clc\ot\clc\ot...\ot\clc}
_{n-times}\raro\underbrace{\clc\ot\clc\ot...\ot\clc}_{n-times}$ is the flip map
between $i$ and $j$-th place and
$m_{ij}:\underbrace{\clc\ot\clc\ot...\ot\clc}_{n-times}
\raro\underbrace{\clc\ot\clc\ot...\ot\clc}_{(n-1)-times}$ is the map obtained by
multiplying $i$ and $j$-th entry. In case we have two copies of an algebra we
shall simply denote by $\sigma$ and $m$ for the flip and multiplication map
respectively.
\subsection{Locally convex $\ast$ algebras and their tensor products}
We begin by recalling from \cite{Takesaki} the tensor product of two $C^{\ast}$
algebras $\clc_{1}$ and $\clc_{2}$ and let us choose the minimal or spatial
tensor products between two $C^{\ast}$ algebras. The corresponding $C^{\ast}$
algebra will be denoted by $\clc_{1}\hat{\ot}\clc_{2}$ throughout this paper.
However we need to consider more general topological spaces and algebras. A
locally convex space is a vector space equipped with a locally convex topology
given by a family of seminorms.  We call a
locally convex space Fr\'echet if the family of seminorms is countable
(hence the space
is metrizable) and is complete with respect to the metric given by the family
of seminorms. There are many ways to equip the algebraic
tensor product of two locally convex spaces with a locally convex topology. Let
$E_{1}$, $E_{2}$ be two
locally convex spaces with family of seminorms $\{||.||_{1,i}\}$ and
$\{||.||_{2,j}\}$ respectively. Then one wants a family $\{||.||_{i,j}\}$ of
seminorms for $E_{1}\ot E_{2}$ such that $||e_{1}\ot
e_{2}||_{i,j}=||e_{1}||_{1,i}||e_{2}||_{2,j}$. The problem is that such a choice
is far from unique and there is a maximal and a minimal choice giving the
projective and injective tensor product respectively. A Fr\'echet locally convex space is called nuclear if its projective
and
injective tensor products with any other Fr\'echet space coincide as a
locally convex space. We do not go into further details
of this topic
here but refer the reader to \cite{TVS} for a comprehensive discussion. 
Furthermore if the space is a $\ast$ algebra then we demand that its $\ast$
algebraic structure is compatible with its locally convex topology i.e.
the involution $\ast$ is continuous and multiplication is jointly continuous with
respect to the topology. Projective and
injective tensor product of two such topological $\ast$ algebras are again
topological $\ast$ algebra. We shall mostly use unital $\ast$
algebras. Henceforth all the topological $\ast$-algebras will be unital unless
otherwise mentioned.

We now
specialize to a particular class of locally convex $\ast$ algebras called smooth $C^*$-normed algebras defined and studied by Blackadar and Cuntz in \cite{smooth}.
 These are $C^*$-normed $\ast$-algebras which are complete w.r.t. the  locally convex topology given by all closable derived seminorms (in the sense of \cite{smooth}). It is proved in \cite{smooth} that such algebras are 
   closed under holomorphic functional calculus, and embedded as a norm-dense $\ast$-subalgebra in a (unique upto isomorphism) $C^*$ algebra. Moreover, 
    any unital $\ast$-homomorphism between such algebras is automatically continuous w.r.t. the corresponding locally convex topologies.
    
    We actually need a slightly smaller subclass of such algebras, to be called `nice algebra' for convenience. This is one of special class of examples of Fr\'echet ($D^{\ast}_{\infty}$)-subalgebra of $C^{\ast}$ algebra treated in \cite{Bhatt}.
\bdfn
 A unital Fr\'echet $\ast$ algebra $\cla$ will be called a `nice' algebra if there is a $C^*$-norm $\| \cdot \|$ on $\cla$ and the 
  the underlying locally convex topology of $\cla$ 
  comes from the family of seminorms $ \| \cdot \|_{\underline{\alpha}} $, with 
$\underline{\alpha}=(i_{1},...,i_{k})$  being any finite multi
index (including the empty index, i.e. $\alpha=\phi$) and where 
$$||x||_{\underline{\alpha}}:=||\delta_{\underline{\alpha}}(x)|| \equiv \| \delta_{i_{1}}\ldots \delta_{i_{k}}(x)\|,$$
$\delta_{\phi}:={\rm id}$ and each $\delta_i$ denotes a $\| \cdot \|$-closable $\ast$-derivation from $\cla$ to itself.
\edfn
 Given two such `nice' algebras
$\cla(\subset\cla_{1})$ and $\clb(\subset\clb_{1})$, where $\cla_1, \clb_1$ denote respectively the $C^*$-completion of $\cla,\clb$ in the corresponding $C^*$-norms,
 we
choose the injective tensor product norm on $\cla\ot \clb$, i.e. we view it as a dense subalgebra of $ \cla_{1}\hat{\ot}\clb_{1}$. $\cla\ot\clb$ has
natural (closable $\ast$) derivations of the forms $\widetilde{\delta}=\delta \ot {\rm id}$ as well as $\widetilde{\eta}={\rm id} \ot \eta$ 
 where $\delta,\eta$ are closable $\ast$-derivations on $\cla$ and $\clb$ respectively. Clearly, $\widetilde{\delta}$ commute with $\widetilde{\eta}$.
 We topologize $\cla\ot \clb$ by the family
of seminorms coming from such derivations, i.e. $\{||.||_{\underline{\alpha}\underline{\beta}}\}$  where
$\| \cdot \|$ is the injective $C^*$-norm and 
$$ ||X||_{\underline{\alpha}\underline{\beta}}=||\widetilde{\delta}_{\underline{\alpha}}\widetilde{\eta}_{\underline{\beta}}(X)||,$$ $\underline{\alpha}
=(i_1,\ldots, i_k)$, $\underline{\beta}=(j_1,\ldots, j_l)$ some multi indices as before and $\delta_i, \eta_j$'s being 
closable $\ast$-derivations on $\cla$ and $\clb$ respectively.

We denote the completion of $\cla\ot\clb$ with respect to this topology by $\cla\hat{\ot}\clb$. It follows from \cite{smooth} that it is indeed a $C^*$-normed smooth algebra and 
it is also Fr\'echet. However, we cannot in general argue that it is nice in our sense as there may be more derivations on the completion which do not `split' as a sum of 
 derivations on the two constituent algebras. 
 Fortunately, all the
locally convex
$\ast$ algebras considered in this paper will be of the form
$C^{\infty}(M)\hat{\ot}\clq$ for some compact manifold $M$ and some unital  $C^{\ast}$
algebra $\clq$ so that they will be `nice' locally convex $\ast$ algebras in
our sense with $C(M)\hat{\ot}\clq$ as the ambient $C^{\ast}$ algebra and the
finitely many canonical derivations coming from the coordinate vector fields
of the compact manifold $M$. We will show in one of the appendices that the above tensor product $\hat{\ot}$ of such 
 algebras will turn out to be nice again.\\
  \indent In case a nice algebra $\clb$  is complete in the $C^*$-norm, i.e. a $C^*$-algebra, all its $\ast$-derivations are norm-bounded. Thus, for any other 
  nice algebra $\cla$, the topology of the tensor product $\cla \hat{\ot} \clb$ is determined by $\widetilde{\delta}$'s only, where $\delta$'s are closable $\ast$-derivations 
   on $\cla$.\\ 
\indent Given  nice $\ast$ algebras $\cle_{1},\cle_{2},\clf_{1},\clf_{2}$ 
 and
$u:\cle_{1}\raro \cle_{2}$, $v:\clf_{1}\raro \clf_{2}$ two unital $\ast$-
homomorphisms which are automatically continuous, the algebraic tensor product map $u\ot v$ can be shown to be 
 continuous with respect to the locally
convex topology discussed above.  We denote the continuous extension again by
$u{\ot}v.$ Using this, we can also define $({\rm id} \ot \omega): \cle_1 \hat{\ot} \cle_2 \raro \cle_1$ if $\cle_2$ is a $C^*$ algebra and  $\omega$ is any state, or more generally,
 a bounded linear functional.

We note the following standard fact without proof which will be crucial in the analysis of smooth actions of compact quantum groups later on.
\bppsn
If $\cla_1$, $\cla_2,$ $\cla_3$ are nice algebras as above and $\Phi: \cla_1
\times \cla_2 \raro \cla_3$ is a bilinear map which is separately continuous in
each of the arguments. Then $\Phi$ extends to a continuous linear map from the
projective tensor product of $\cla_1$ with $\cla_2$ to $\cla_3$. If furthermore,
$\cla_1$ is nuclear, $\Phi$ extends to a continuous map from $\cla_1 \hat{\ot}
\cla_2$ to $\cla_3$. 

As a special case, suppose that $\cla_1,\cla_2$ are subalgebras of a nice algebra $\cla$ and also that $\cla_1$ is isomorphic as a Fr\'echet space to some nuclear space. 
Then the multiplication map, say $m$, of $\cla$ extends to a continuous map from $\cla_1 \hat{\ot} \cla_2$ to $\cla$.
\eppsn
   
\subsection{Exterior tensor product of Hilbert bimodules}
Let $\cle_{1}$ and $\cle_{2}$ be two pre Hilbert bimodules over two locally
convex $*$ algebras $\clc_{1}$ and $\clc_{2}$ which are subalgebras of nice algebras respectively. We denote the
algebra
valued inner product for the pre Hilbert bimodules by $<<,>>$. When the bimodule is
a
pre Hilbert space, we denote the corresponding scalar valued inner product by
$<,>$. Then $\cle_{1}\ot \cle_{2}$ has an obvious $\clc_{1}\ot \clc_{2}$
bimodule
structure, given by $(a\ot b)(e_{1}\ot e_{2})(a^{\prime}\ot b^{\prime})=ae_{1}a^{\prime}\ot
be_{2}b^{\prime}$ for $a,a^{\prime}\in \clc_{1}, b,b^{\prime}\in \clc_{2}$ and $e_{1}\in
\cle_{1},
e_{2}\in \cle_{2}$. Also define $\clc_{1}\ot \clc_{2}$ valued inner
product by $<<e_{1}\ot e_{2},f_{1}\ot f_{2}>>=<<e_{1},f_{1}>>\ot
<<e_{2},f_{2}>>$ for $e_{1}, f_{1}\in \cle_{1}$ and $e_{2},f_{2}\in \cle_{2}$. In case $\clc_{1}$,$\clc_{2}$ are complete, i.e. nice algebras and $\cle_{1},\cle_{2}$ are Hilbert bimodules, 
we denote the completed module by $\cle_{1}\bar{\ot}\cle_{2}$. In fact
$\cle_{1}\bar{\ot}\cle_{2}$ is an
$\clc_{1}\hat{\ot}\clc_{2}$ bimodule. This is called the exterior tensor product
of two
bimodules. In particular if one of the bimodule is a Hilbert space $\clh$
(bimodule over $\mathbb{C}$) and the other is a $C^{*}$ algebra $\clq$
(bimodule over itself), then the exterior tensor product gives the usual
Hilbert $\clq$ module $\clh\bar{\ot}\clq$. When
$\clh=\mathbb{C}^{N}$, we have a natural identification of an element
$T=((T_{ij}))\in M_{N}(\clq)$ with the right $\clq$ linear map of
$\mathbb{C}^{N}\ot \clq$ given by
$$e_{i}\mapsto \sum e_{j}\ot T_{ji},$$
where $\{e_{i}\}_{i=1,...,N}$ is a basis for $\mathbb{C}^{N}$. We shall need tensor product of maps, which are `isometric' in some sense. Let
$T_{i}:\cle_{i}\raro\clf_{i}$, $i=1,2$ be two $\mathbb{C}$-linear maps and
$\cle_{i}$, $\clf_{i}$ be Hilbert bimodules over $\clc_{i}$, $\cld_{i}$
($i=1,2$) respectively. Moreover, suppose that
$<<T_{i}(\xi_{i}),T_{i}(\eta_{i})>>=\alpha_{i}<<\xi_{i},\eta_{i}>>$,
$\xi_{i},\eta_{i}\in \cle_{i}$ where $\alpha_{i}:\clc_{i}\raro\cld_{i}$ are
$\ast$-homomorphisms. Then it is easy to show that the algebraic tensor product
$T:=T_{1}\ot_{alg}T_{2}$ also satisfies
$<<T(\xi),T(\eta)>>=(\alpha_{1}\ot\alpha_{2})(<<\xi,\eta>>)$ and hence extends
to a well defined continuous map from $\cle_{1}\bar{\ot}\cle_{2}$ to
$\clf_{1}\bar{\ot}\clf_{2}$ again to be denoted by $T_{1}\ot T_{2}$. 
\subsection{Compact quantum groups, their representations and actions}
\subsubsection{Definition and representations of compact quantum groups}
A compact quantum group (CQG for short) is a  unital $C^{\ast}$ algebra $\clq$ with a
coassociative coproduct 
(see \cite{Van}, \cite{Pseudogroup}) $\Delta$ from $\clq$ to $\clq \hat{\ot} \clq$  
such that each of the linear spans of $\Delta(\clq)(\clq\ot 1)$ and that
of $\Delta(\clq)(1\ot \clq)$ is norm-dense in $\clq \hat{\ot} \clq$. 
From this condition, one can obtain a canonical dense unital $\ast$-subalgebra
$\clq_0$ of $\clq$ on which linear maps $\kappa$ and 
$\epsilon$ (called the antipode and the counit respectively) are defined making the above subalgebra a Hopf $\ast$ algebra. In fact, this is  the algebra generated by the `matrix coefficients' of
the (finite dimensional) irreducible non degenerate representations (to be
defined 
 shortly) of the CQG. The antipode is an anti-homomorphism and also satisfies $\kappa(a^{\ast})=(\kappa^{-1}(a))^{\ast}$ for $a \in \clq_{0}$.
 
 It is known that there is a unique state $h$ on a CQG $\clq$ (called the Haar
state) which is bi invariant in the sense that $({\rm id} \ot h)\circ
\Delta(a)=(h \ot {\rm id}) \circ \Delta(a)=h(a)1$ for all $a$. The Haar state
need not be faithful in general, though it is always faithful on $\clq_0$ at
least. Given the Hopf $\ast$-algebra $\clq_{0}$, there can be several CQG's
which have this $\ast$-algebra as the Hopf $\ast$-algebra generated by the
matrix elements of finite dimensional representations. We need two of such CQG's: the reduced and the universal one. By definition, the reduced 
 CQG $\clq_r$ is the image of $\clq$ in the GNS representation of $h$, i.e. $\clq_r=\pi_r(\clq)$, $\pi_r: \clq \raro \clb(L^2(h))$ is the GNS representation. 
 
 There also exists a
largest such CQG $\clq^{u}$, called the universal CQG corresponding to
$\clq_{0}$. It is obtained as the universal enveloping $C^{\ast}$ algebra of
$\clq_{0}$. We also say that a CQG $\clq$ is universal if $\clq=\clq^{u}$. Given two CQG's $(\clq_{1},\Delta_{1})$ and $(\clq_{2},\Delta_{2})$, a $\ast$ homomorphism $\pi\clq_{1}\raro\clq_{2}$ is said to be a CQG morphism if $(\pi\ot\pi)\circ\Delta_{1}=\Delta_{2}\circ\pi$ on $\clq_{1}$. In case $\pi$ is surjective, $\clq_{2}$ is said to be a quantum subgroup of $\clq_{1}$ (or a Woronowicz subalgebra) and denoted by $\clq_{2}\leq\clq_{1}$.\\
%
\indent Let $\clh$ be a Hilbert space. Consider the multiplier algebra
$\clm(\clb_{0}(\clh)\hat{\ot} \clq)$. This algebra has two natural embeddings
into $\clm(\clb_{0}(\clh)\hat{\ot} \clq\hat{\ot} \clq)$. The first
one is
obtained by extending the map $x\mapsto x\ot 1$. The second one is obtained by
composing this map with the flip on the last two factors. We will write $w^{12}$
and $w^{13}$ for the images of an element $w\in \clm(\clb_{0}(\clh)\hat{\ot}
\clq)$
by these two maps respectively. Note that if $\clh$ is finite dimensional then
$\clm(\clb_{0}(\clh)\hat{\ot} \clq)$ is isomorphic to $\clb(\clh)\ot \clq$
(we do not need any topological completion).
\bdfn
Let $(\clq,\Delta)$ be a CQG. A unitary representation of $\clq$ on a Hilbert
space $\clh$ is a $\mathbb{C}$-linear map $U$ from $\clh$ to the Hilbert module
$\clh\bar{\ot}
\clq$ such that \\
1. $<<U(\xi),U(\eta)>>=<\xi,\eta>1_{\clq}$, where $\xi,\eta\in \clh$.\\
2. $(U\ot {\rm id})U=({\rm id}\ot \Delta)U.$
 \edfn
\indent Given such a unitary representation we have a unitary 
element $\widetilde{U}$ belonging to
$\clm(\clb_{0}(\clh)\hat{\ot}
\cls)$ given by $\widetilde{U}(\xi\ot b)=U(\xi)b,(\xi\in \clh, b\in \cls)$
satisfying 
$({\rm id}\ot \Delta)\widetilde{U}= \widetilde{U}^{12}\widetilde{U}^{13}$, where $\clm(\clc)$ denotes the multiplier algebra of a $C^{\ast}$ algebra $\clc$.
\bdfn
A closed subspace $\clh_{1}$ of $\clh$ is said to be invariant if
$U(\clh_{1})\subset \clh_{1}\bar{\ot}\clq$. A unitary representation $U$ of a CQG
is said to be irreducible if there is no proper
invariant subspace.
\edfn
It is a well known fact that every irreducible unitary representation is finite
dimensional.\\
\indent We denote by $Rep(\clq)$ the set of inequivalent irreducible unitary
representations
of $\clq$. For $\pi\in Rep(\clq)$, let $d_{\pi}$ and $\{t^{\pi}_{jk}:
j,k=1,...,d_{\pi} \}$ be the dimension and matrix coefficients of the
corresponding finite dimensional representation respectively. Then for each
$\pi\in Rep(\clq)$, we have a unique $d_{\pi}\times d_{\pi}$ complex
matrix $F_{\pi}$ such that\\
(1) $F_{\pi}$ is positive and invertible with $Tr(F_{\pi})=Tr
(F_{\pi}^{-1})=M_{\pi}>0$(say).\\
(2) $h(t_{ij}^{\pi}t_{kl}^{\pi^{\ast}})=
\frac{1}{M_{\pi}}\delta_{ik}F_{\pi}(j,l).$\\ 
\indent Corresponding to $\pi\in Rep(\clq)$, let $\rho^{\pi}_{sm}$ be the linear
functional on $\clq$ given by $\rho^{\pi}_{sm}(x)=h(x^{\pi}_{sm}x),
s,m=1,...,d_{\pi}$ for $x\in \clq$ where $x^{\pi}_{sm}=(M_{\pi})t^{\pi
*}_{km}(F_{\pi})_{ks}$. Also let
$\rho^{\pi}=\sum_{s=1}^{d_{\pi}}\rho^{\pi}_{ss}$.\\
\indent We say a map $\Gamma:\clk\raro \clk\ot \clq_{0}$ (where $\clk$ is a vector space apriori without any topology) is an algebraic representation of the CQG $\clq$ if $\Gamma$ is algebraic, $(\Gamma\ot {\rm id})\Gamma=({\rm id}\ot
 \Delta)\Gamma$ and ${\rm Sp} \ \Gamma(\clk)\clq_{0}=\clk\ot\clq_{0}$.  
\subsubsection{Actions of compact quantum groups}
Let $\clc$ be a nice  unital Fr\'echet $\ast$-algebra (in the sense discussed
in Subsection 3.1) and $\clq$ be a
compact
quantum group. 
\bdfn
\label{CQG_top_action}
A $\mathbb{C}$ linear map $\alpha:\clc\raro \clc\hat{\ot}\clq$ is
said to be a topological action of $\clq$ on $\clc$ if\\
1. $\alpha$ is a continuous $\ast$ algebra homomorphism.\\
2. $(\alpha\ot {\rm id})\alpha=({\rm id}\ot \Delta)\alpha$ (co-associativity).\\
3. $Sp \ \alpha(\clc)(1\ot \clq)$ is dense in $\clc\hat{\ot}\clq$ in the
corresponding Fr\'echet topology.
\edfn
Note that if the Fr\'echet algebra is a $C^{\ast}$ algebra, then the
definition of a topological action coincides with the usual $C^{\ast}$ action
of a compact quantum group.
\bdfn
A topological action $\alpha$ is said to be faithful if the $\ast$-subalgebra
of $\clq$ generated by the elements of the form $(\omega\ot {\rm id})\alpha$, where
$\omega$ is a continuous linear functional on $\clc$, is dense in $\clq$.
\edfn
 
\bdfn
Let $(H,\Delta,\epsilon,\kappa)$ be a Hopf $\ast$ algebra and $A$ be a $\ast$ algebra. A unital $\ast$ algebra homomorphism $\alpha:A\raro A\ot H$ is said to be a Hopf $\ast$ algebraic (co)action of $H$ on $A$ if\\
(i) $(\alpha\ot{\rm id})\alpha=({\rm id}\ot \Delta)\alpha.$\\
(ii) $({\rm id}\ot \epsilon)\alpha={\rm id}.$
\edfn
For a Hopf $\ast$ algebraic (co)action as above, one can prove that ${\rm Sp} \ \{\alpha(A)(1\ot H)\}=A\ot H$. 
\brmrk
For a topological action $\alpha$ of the CQG $\clq$ on the nice algebra $\clc$, we say that it is algebraic over a $\ast$-subalgebra $\clc_{0}\subset\clc$ if $\alpha|_{\clc_{0}}:\clc_{0}\raro\clc_{0}\ot\clq_{0}$ is a Hopf $\ast$-algebraic action. 
\ermrk
\indent Let $X$ be a compact space. Then we can consider the $C^{\ast}$ action
of $\clq$ on $C(X)$. We say $\clq$ acts topologically on a compact space $X$ if
there is a $C^{\ast}$ action of $\clq$ on $C(X)$. We have the following(\cite
{Huichi_huang}):
\bppsn
\label{tracial}
If a CQG $\clq$ acts topologically and faithfully on $X$, where $X$ is any
compact space, then
the corresponding reduced CQG $\clq_{r}$ (which has a faithful Haar state) must
be
a Kac algebra. In particular the Haar state of $\clq_{r}$, and hence of $\clq$
is
tracial. Moreover, the antipode $\kappa$ is defined and norm-bounded on $\clq_r$.
\eppsn   
 \indent Now recall from Subsection 4.1 the linear functional $\rho^{\pi}$ on a
 compact quantum group $\clq$. Given a topological action $\alpha$ of $\clq$
 on a `nice' unital $\ast$ algebra $\clc$, we can define a projection $P_{\pi}:\clc\raro \clc$ called spectral projection corresponding to $\pi\in Rep(\clq)$
  by $P_{\pi}:=({\rm id}\ot \rho^{\pi})\alpha$. Note that
 $({\rm id}\ot \phi)\alpha(\clc)\subset \clc$ for all bounded linear functionals
 $\phi$ on $\clq$. We call $\clc_{\pi}:={\rm Im} \ P_{\pi}$ the spectral subspace corresponding to $\pi$. The subspace spanned by $\{\clc_{\pi}, \pi\in Rep(\clq)\}$ is actually a unital $\ast$-subalgebra called the spectral subalgebra for the action. Let $\clc_{0}:=ker(\alpha)\oplus \ Sp
\{\clc_{\pi}:\pi\in\hat{\clq}\} $. So in particular if $ker(\alpha)=\{0\}$,
then the spectral subspace coincides with $\clc_{0}$. Along the lines of
\cite{Podles} and Proposition 2.2 of \cite{Act_com} we have 
\bppsn
\label{maximal}
 (i) $\clc_{0}$ is a unital $\ast$ subalgebra over which $\alpha$ is algebraic.\\
(ii) $\clc_{0}$ is dense in $\clc$ in the Fr\'echet topology.\\
(iii) $\clc_{0}$ is maximal among the subspaces $V\subset\clc$ with the property that $\alpha(V)\subset V\ot\clq_{0}$.
\eppsn
We end this subsection with a discussion on an analogue of an action on a von
Neumann algebra given by conjugating unitary representation. Let
$\clm\subset\clb(\clh)$ be a von Neumann algebra and $U$ be a unitary
representation of a CQG $\clq$ on $\clh$. The map ${\rm ad}_{\widetilde{U}}$ is a normal,
injective $\ast$-homomorphism on $\clb(\clh)$. We say ${\rm ad}_{\widetilde{U}}$ leaves
$\clm$ invariant if $({\rm id}\ot\phi)({\rm ad}_{\widetilde{U}}(x))\in \clm$ for all bounded
linear functionals $\phi$ on $\clq$. In that case we can consider the spectral
projections $P_{\pi}=({\rm id}\ot \rho^{\pi})\circ {\rm ad}_{\widetilde{U}}$ and define
$\clm_{\pi}:= P_{\pi}(\clm)$. As before, we have the spectral subalgebra $\clm_{0}= \ Sp \
\{\clm_{\pi}, \pi \in \hat{\clq} \}$ and as ${\rm ad}_{\widetilde{U}}$ is one to one, $\clm_{0}$ is maximal among subspaces $V\subset\clm$ for which ${\rm ad}_{\widetilde{U}}(V)\subset V\ot\clq_{0}$. Furthermore, it can be proved that $\clm_{0}$ is SOT
dense in $\clm$. 
 We end this subsection with a discussion on the unitary implementability of an action. 
 \bdfn
 We call an action $\alpha$ of a CQG $\clq$ on a unital $C^*$ algebra $\clc$ to be implemented by a unitary representation $U$ of $\clq$ in $\clh$, say, if
  there is a faithful representation $\pi : \clc \raro \clb(\clh)$ such that $\widetilde{U} (\pi(x) \ot 1)\widetilde{U}^*=(\pi \ot {\rm id})(\alpha(x))$ for all $x \in \clc$.
  \edfn
  It is clear that if an action is implemented by a unitary representation then it is one-to-one. In fact, as $({\rm id} \ot \pi_r)(U)$ gives a unitary representation of $\clq_r$ in $\clh$
   and the `reduced action' $\alpha_r:=({\rm id} \ot \pi_r)\circ \alpha$ of $\clq_r$ is also implemented by a unitary representation, it follows that even 
    $\alpha_r$ is one-to-one. We see below that this is actually equivalent to implementability by unitary representation:\\
    \blmma
    \label{unitary_impl}
    Given an action $\alpha$ of $\clq$ on a unital separable $C^*$ algebra $\clc$ the following are equivalent:\\
    (a) There is a faithful positive functional $\phi$ on $\clc$ which is invariant w.r.t. $\alpha$, i.e. $(\phi \ot {\rm id})(\alpha(x))=\phi(x)1_\clq$ for all $x \in \clc$.\\
    (b) The action is implemented by some unitary representation.\\
    (c) The reduced action $\alpha_r$ of $\clq_r$ is injective.\\
\elmma
{\it Proof:}\\
If (a) holds, we consider $\clh$ to be the GNS space of the faithful positive functional $\phi$. The GNS representation $\pi$ is faithful, and the linear map $U$ defined by 
 $U(x):=\alpha(x)$ from $\clc \subset \clh=L^2(\clc,\phi)$ to $\clh \bar{\ot} \clq$ is an isometry by the invariance of $\phi$. Thus $U$ extends to $\clh$ and it is easy to check that 
  it gives a unitary representation which implements $\alpha$.\\
  
  We have already argued $(b) \Rightarrow (c)$, and finally, if (c) holds, we choose any faithful state say $\tau$ on the separable $C^{\ast}$ algebra $\clc$ and take $\phi(x)=(\tau \ot h)(\alpha_r(x))$,
   which is faithful as $h$ is faithful on $\clq_r$. It can easily be verified that $\phi$ is $\alpha$-invariant on the dense subalgebra $\clc_0$ mentioned before, and hence on the whole of 
    $\clc$.\qed
 
\subsection{\bf Representation of CQG on a Hilbert bimodule over a
nice topological $\ast$-algebra}
We now generalize the notion of unitary representation on Hilbert
spaces in
another direction, namely on Hilbert bimodules over nice, unital topological
$\ast$-algebras. Let $\cle$ be a Hilbert $\clc-\cld$ bimodule over topological
$\ast$-algebras $\clc$ and $\cld$. Also let $\clq$ be a compact quantum group. If
we consider
$\clq$ as a bimodule over itself, then we can form the exterior tensor product
$\cle\bar{\ot}\clq$ which is a
$\clc\hat{\ot} \clq-\cld\hat{\ot}\clq$ bimodule. Also let
$\alpha_{\clc}:\clc\raro \clc\hat{\ot} \clq$ and $\alpha_{\cld}:\cld\raro
\cld\hat{\ot} \clq$ be topological actions on $\clc$ and $\cld$ of $\clq$ in the
sense discussed earlier. Using $\alpha$ we can
give $\cle\bar{\ot}\clq$ a $\clc-\cld$ bimodule structure given by
$a.\eta. a^{'}= \alpha_{\clc}(a)\eta \alpha_{\cld}(a^{\prime})$, for $\eta\in
\cle\bar{\ot}\clq$ and $a\in\clc,a^{\prime}\in \cld$ (but without any $\cld$ valued
inner
product).
\bdfn
A $\mathbb{C}$-linear map $\Gamma:\cle\raro \cle\bar{\ot}\clq$
is
said to be an $\alpha_{\cld}$ equivariant unitary representation of $\clq$ on
$\cle$ if\\
1. $\Gamma(\xi d)=\Gamma(\xi)\alpha_{\cld}(d)$ and
$\Gamma(c\xi)=\alpha_{\clc}(c)\Gamma(\xi)$ for $c\in\clc,d\in\cld$.\\
2. $<<\Gamma(\xi),\Gamma(\xi^{\prime})>>=\alpha_{\cld}(<<\xi,\xi^{\prime}>>)$, for
$\xi,\xi^{\prime}\in \cle$.\\
3. $(\Gamma\ot {\rm id})\Gamma=({\rm id}\ot \Delta)\Gamma$ (co associativity)\\
4. $\overline{Sp} \ \Gamma(\cle)(1\ot \clq)=\cle\bar{\ot}\clq$ (non
degeneracy).
\edfn
In the above definition note that condition (2) allows one to define
$(\Gamma\ot {\rm id})$. If $\clc=\cld$ and $\alpha_{\clc}=\alpha_{\cld}$ we simply call $\Gamma$ $\alpha$-equivariant. If $\alpha_{\clc}$ and $\alpha_{\cld}$ are understood from the context, we may call $\Gamma$ just equivariant. Given an $\alpha$ equivariant representation $\Gamma$ of
$\clq$ on a Hilbert
bimodule $\cle$, proceeding as in Subsection 4.2, we define $P_{\pi}:=({\rm id}\ot \rho^{\pi})\circ \Gamma$, $\cle_{\pi}:={\rm Im}P_{\pi}$ for $\pi\in Rep(\clq)$ and $\cle_{0}:=Sp \ \{\cle_{\pi};\pi\in \ Rep(\clq)\}\oplus ker(\Gamma)$. In case $\Gamma$ is one-one
which is equivalent to $\alpha_{\cld}$ being one-one, $\cle_{0}$ coincides with $Sp \ \{\cle_{\pi};\pi\in \ Rep(\clq)\}$. Again proceeding along the lines of
\cite{Podles} and \cite{Act_com}, we can prove the following analogue of Proposition \ref{maximal}:
\bppsn
\label{sp_dec}
 1. $({\rm id}\ot \phi)\Gamma(\cle)\subset \cle$, for all bounded linear functional
$\phi\in \clq^{\ast}$.\\
 2. $P_{\pi}^{2}=P_{\pi}.$\\
 3. $\Gamma$ is algebraic over $\cle_{0}$ i.e. $\Gamma(\cle_{0})\subset
\cle_{0}\ot
 \clq_{0}$ and $({\rm id}\ot \epsilon)\Gamma= {\rm id}$ on $\cle_{0}$.\\ 
 4. $\cle_{0}$ is dense in $\cle$.\\
 5. $\cle_{0}$ is maximal among the subspaces $V\subset\cle$ such that $\Gamma(V)\subset V\ot\clq$.
\eppsn

\subsection{Locally convex $\ast$ algebras and Hilbert bimodules coming from
classical geometry}
 
\subsubsection{$C^{\infty}(M)$ as a nice algebra}
All the notations are as in the Subsection 3.1 and throughout the section, we denote a smooth $n$-
dimensional compact manifold
 possibly with
 boundary by $M$. We
 denote the algebra of real (complex respectively) valued smooth
 functions on $M$ by $C^{\infty}(M)_{\mathbb{R}}$ $(C^{\infty}(M) \ {\rm respectively})$. Clearly
 $C^{\infty}(M)$ is the complexification of $C^{\infty}(M)_{\mathbb{R}}$ . We
 also equip it with a locally convex topology : we say a sequence $f_{n}\in
 C^{\infty}(M)$ converges to an $f\in C^{\infty}(M)$ if for every compact set
$K$
 within a single
 coordinate neighborhood ($M$ being compact, has finitely many such
 neighborhoods) and a multi-index $\alpha$, $\partial^{\alpha}f_{n}
 \raro  \partial^{\alpha}f$ uniformly over $K$. Here, for a coordinate neighborhood $(x_1, \ldots, x_n)$ and a multi-index
  $\alpha=(i_1,\ldots, i_k)$ with $i_j \in \{ 1, \ldots, n\}$,  $\partial^{\alpha}$ denotes $\frac{\partial}{\partial x_{i_1}} \ldots \frac{\partial}{\partial x_{i_k}}$.
  Equivalently, let $U_{1},
 U_{2},...,U_{l}$ be a finite cover of $M$. Then it is a locally
 convex
 topology described by a countable family of seminorms given by: 
 $$p_{i}^{K,\alpha}(f)=\sup_{x\in K}|\partial^{\alpha}f(x)|,$$ where $K$ is a
 compact
 set within $U_{i}$, $\alpha$ is any multi index, $i=1,2,....l$.
$C^{\infty}(M)$ is complete with respect to this topology (Example 1.46 of
\cite{Rudin} with obvious modifications) and hence this makes $C^{\infty}(M)$ a
locally convex Fr\'echet $\ast$ algebra with obvious $\ast$ structure.  In fact, it is nice algebra because any closable $\ast$-derivation on $C^\infty(M)$ comes from 
 a smooth vector field, i.e. locally a $C^\infty(M)$-linear combination of partial derivatives.
 
 Actually, 
by choosing a finite $C^{\infty}$ partition of
unity on the compact manifold $M$, we can obtain finite set
$\{\delta_{1},...,\delta_{N}\}$ for some $N\geq n$ of globally defined vector
fields on
$M$ which is complete in the sense that $\{\delta_{1}(m),...,\delta_{N}(m)\}$
spans $T_{m}(M)$ for all
$m$ (need not be a basis). It follows that for describing the locally convex topology on
$C^{\infty}(M)$ it is enough to consider the seminorms $\| \delta_{i_1}\ldots \delta_{i_k} (\cdot)\|$, with $k \geq 0$, $i_j \in \{ 1, \ldots, N \}$.

It is also known that (see Example 6.2 of \cite{nuclear}
\bppsn
With the above 
topology  $C^{\infty}(M)$ is  nuclear as a locally convex space.\eppsn

 \indent  Let $E$ be any locally convex space. Then we can define the space of
$E$ valued
 smooth functions on a compact manifold $M$. Take a centered coordinate chart
 $(U,\psi)$ around a point $x\in M$. Then an $E$ valued function $f$ on $M$
 is said to be smooth at $x$ if $f\circ \psi^{-1}$ is smooth $E$ valued
 function at $0\in \mathbb{R}^{n}$ in the sense of $\cite{TVS}$(Definition
 40.1). We denote the space of $E$ valued smooth functions on $M$ 
 by $C^{\infty}(M,E)$. We can give a locally convex topology on
 $C^{\infty}(M,E)$ by the family of seminorms given by $p_{i}^{K,\alpha}(f):=
 sup_{x\in K}||\partial^{\alpha}f(x)||$, where $i,K,\alpha$ are as before. Then
 we have the following
 \bppsn
 1. If $E$ is complete, then so is $C^{\infty}(M,E)$.\\
  2. Suppose $E$ is a complete locally convex space. Then we have
   $C^{\infty}(M)\hat{\ot} E\cong C^{\infty}(M,E).$\\
  3. Let $M$ and $N$ be two smooth compact manifolds with boundary. Then
  $C^{\infty}(M)\hat{\ot} C^{\infty}(N)\cong C^{\infty}(M,C^{\infty}(N))\cong
  C^{\infty}(M\times N
  )$ and contains $C^{\infty}(M)\ot C^{\infty}(N)$ as Fr\'echet dense
subalgebra.\\
  4. Let $\clq$ be a $C^{\ast}$ algebra. Then 
$C^{\infty}(M)\hat
 {\ot}\clq\cong C^{\infty}(M,\clq)$ as Fre\'chet  algebras.
  
 \eppsn
 For the proof of the first statement see 44.1 of \cite{TVS}. The second
statement also follows from 44.1 of \cite{TVS} and the fact that
$C^{\infty}(M)$ is nuclear. The third and fourth
 statements follow from the second statement (replacing $E$ suitably).\\
\indent Although nuclearity ensures that $C^\infty(M) \hat{\ot} \clq$ is unambiguously defined and coincides with 
 the tensor product for nice algebras defined earlier, it remains to see whether it is itself a nice algebra. Indeed, this is true and we'll give a proof in Appendix.
In fact, we'll have more  (see the Appendix):
\blmma
For any two $C^{\ast}$ algebras $\clq, \clq^\prime$, 
$C^{\infty}(M,\clq)$ is a nice algebra and the nice algebra tensor product $C^\infty(M, \clq) \hat{\ot}\clq^{\prime}$ is $\ast$-isomorphic with 
$C^{\infty}(M)\hat{\ot}(\clq\hat{\ot}\clq^{\prime})\cong
C^{\infty}(M,\clq\hat{\ot}\clq^{\prime})$
\elmma

\subsubsection{Hilbert bimodule of one-forms on a Riemannian manifold}
 Let
$\Lambda^{k}(C^{\infty}(M))$ be the space of smooth $k$ forms on the manifold
$M$, with the natural locally convex
topology induced by the topology of $C^{\infty}(M)$ given by a
family of seminorms $q_{i}^{K,\alpha}(\omega)={\rm sup_{x\in K, 1\leq j\leq n}}|\partial^{\alpha}f_{j}(x)|$, $K \subset U_i$, where $K,\alpha,\{U_{1},...,U_{l}\}$ are
as before and $\omega|_{U_{i}}=\sum_{j=1}^{n}f_{j}dx_{j}$. It is clear from the definition that the differential map
$d:C^{\infty}(M)\raro \Omega^{1}(C^{\infty}(M))$ is Fr\'echet continuous.\\
\blmma
 \label{Frechet_den}
 Let $\cla$ be a Fr\'echet dense subalgebra of $C^{\infty}(M)$. Then
$\Lambda^{1}(\cla):=Sp \ \{fdg:f,g\in\cla\}$ is dense in
$\Lambda^{1}(C^{\infty}(M))$.
 \elmma
 {\it Proof:}\\
 It is enough to approximate $fdg$ where $f,g\in C^{\infty}(M)$ by elements of
$\Omega^{1}(\cla)$. By Fr\'echet density of $\cla$ in $C^{\infty}(M)$ we can
choose sequences $f_{m},g_{m}\in \cla$ such that $f_{m}\raro f$ and
$g_{m}\raro g$ in the Fr\'echet topology, hence by the continuity of $d$ and
the $C^{\infty}(M)$ module multiplication of $\Lambda^{1}(C^{\infty}(M))$, we
see that $f_{m}dg_{m}\raro fdg$ in $\Lambda^{1}(C^{\infty}(M))$.\qed\\

Let $\Omega^{k}(C^{\infty}(M))_{u}$ be the space of universal k-forms on the
manifold $M$ and $\delta$ be the derivation for the universal algebra of forms
for $C^{\infty}(M)$ i.e $\delta: \Omega^{k}(C^{\infty}(M))_{u}\raro
\Omega^{k+1}(C^{\infty}(M))_{u}$(see \cite{Landi} for further details).\\

 By the universal property $\exists$ a surjective bimodule morphism
$\pi\equiv\pi_{(1)}:\Omega^{1}(C^{\infty}(M))_{u}\raro
\Lambda^{1}(C^{\infty}(M))$, such that $\pi(\delta g)=dg.$\\
 $\Omega^{1}(C^{\infty}(M))_{u}$ has a $C^{\infty}(M)$ bimodule structure:\\
 $$f(\sum_{i=1}^{n}g_{i}\delta h_{i}) = \sum_{i=1}^{n}fg_{i}\delta h_{i}$$
 $$(\sum_{i=1}^{n}g_{i}\delta h_{i})f = \sum_{i=1}^{n} (g_{i}\delta (h_{i}f)-g_{i}h_{i}\delta f)$$
 As $M$ is compact, there is a Riemannian structure. Using the Riemannian
structure on $M$ we can equip $\Omega^{1}(C^{\infty}(M))$ with a $C^{\infty}(M)$
valued inner product
$<<\sum_{i=1}^{n}f_{i}dg_{i},\sum_{i=1}^{n}f_{i}^{'}dg_{i}^{'}>> \in C^{\infty}
(M)$ by the following prescription:\\
 for $x\in M$ choose a coordinate neighborhood $(U,x_{1},x_{2},....,x_{n})$
around $x$ such that $dx_{1},dx_{2},...,dx_{n}$ is an orthonormal basis for
$T_{x}^*M$. Note that the topology does not depend upon any particular choice of
the Riemannian metric. Then
 \begin{eqnarray}
<<\sum_{i=1}^{n}f_{i}dg_{i},\sum_{i=1}^{n}f_{i}^{\prime}dg_{i}^{\prime}>>(x) =
(\sum_{i,j,k,l}\bar {f_{i}}f_{j}^{\prime}(\frac{\bar {\partial g_{i}}}{\partial
x_{k}}\frac{{\partial g_{j}^{\prime}}}{\partial x_{l}}))(x). 
 \end{eqnarray}
 We see that a sequence
$\omega_{n}\raro \omega$ in $\Lambda^{1}(C^{\infty}(M))$ if $<<\omega_{n}-
\omega,\omega_{n}- \omega>> \raro 0$ in Fr\'echet topology of $C^{\infty}(M). $
With this $\Lambda^{1}(C^{\infty}(M))$ becomes a Hilbert module.\\
\section{Smooth action of a CQG on a manifold}
In this section we consider a compact manifold $M$
possibly with boundary, but not necessarily orientable and discuss a notion of
smoothness of a CQG action $\alpha$ of a CQG $\clq$. Moreover we prove that smoothness automatically implies injectivity of 
the reduced action, hence unitary implementability. This will follow from injectivity of the action on $C^\infty(M)$ which we will prove.
\subsection{Definition of a smooth action}
\bdfn
A topological action (in the sense of Definition \ref{CQG_top_action}) of $\clq$ on the
Fr\'echet algebra $C^{\infty}(M)$ is called the smooth action of $\clq$ on the
manifold $M$. In case $M$ has a boundary, we also assume that the closed ideal $\{ f \in C^\infty(M):~f|_{\partial M}=0 \}$ is invariant by the action.
\edfn
\blmma
A smooth action $\alpha$ of $\clq$ on $M$ extends to a $C^{\ast}$ action on $C(M)$ which is denoted by $\alpha$ again. 
\elmma
{\it Proof}: \\
It follows from generalities about smooth $C^*$-normed algebras proved in \cite{smooth} (Proposition 6.8 there), which uses the fact that  $C^{\infty}(M)$ is stable under taking square
roots of positive invertible elements. \qed
\blmma
\label{smooth_ns}
Given a $C^{\ast}$ action $\alpha:C(M)\raro C(M)\hat{\ot}\clq$, 
$\alpha(C^{\infty}(M))\subset C^{\infty}(M,\clq)$ if and only if
$({\rm id}\ot\phi)(\alpha(C^{\infty}(M)))\subset C^{\infty}(M)$ for all bounded linear
functionals $\phi$ on $\clq$.
\elmma{\it Proof}:\\
Only if part:\\
See discussion in Subsection 4.2.\\
If part: \\
follows from Corollary \ref{one} of Appendix (Section \ref{appendix}).\qed

\bthm
\label{suff_smooth}
Suppose we are given a $C^{\ast}$ action $\alpha$ of $\clq$ on $M$. Then
following are equivalent:\\
1) $\alpha(C^{\infty}(M))\subset C^{\infty}(M,\clq)$ and $Sp \
\alpha(C^{\infty}(M))(1\ot \clq)$ is Fr\'echet dense in $C^{\infty}(M,\clq).$\\
2) $\alpha$ is smooth.\\
3) $({\rm id}\ot \phi)\alpha(C^{\infty}(M))\subset C^{\infty}(M)$ for every state $\phi$
on $\clq$, and there is a Fr\'echet dense subalgebra $\cla$ of $C^{\infty}(M)$
over
which $\alpha$ is algebraic.\\
\ethm
{\it Proof:}\\
(1)$\Rightarrow$ (2): Observe that it is enough to show that $\alpha$ is Fr\'echet continuous. But this follows from the Closed Graph Theorem, as 
 Fr\'echet topology is stronger than the norm topology.\\
(2)$\Rightarrow$ (3): Follows from the Proposition \ref{maximal}.\\
(3)$\Rightarrow$ (1): From Lemma \ref{smooth_ns}, it follows that $\alpha(C^{\infty}(M))\subset
C^{\infty}(M,\clq)$. The density condition follows from densities of $\cla$
 and $\cla\ot \clq_{0}$ in $C^{\infty}(M)$ and $C^{\infty}(M,\clq)$
respectively.\qed\\
\subsection{Injectivity of the smooth action}
In this subsection we show that for a smooth action $\alpha:C^{\infty}(M)\raro C^{\infty}(M)\hat{\ot}\clq$, the corresponding reduced $C^*$ action on $C(M)$ is injective. Note that for a general $C^*$ action this is not true. For a CQG $\clq$ where $\clq$ is a non amenable $C^*$ algebra, the coproduct gives an action of the reduced CQG which is not injective (see \cite{Act_com}). We begin by proving an interesting fact which will be used later. For that recall that $C^\infty(M)$ is a nuclear locally convex space and hence so is any quotient by closed ideals. 
\blmma
\label{smooth_bdd_counit}

If $\clq$ has a faithful smooth action on $C^\infty(M)$, where $M$ is compact manifold, then for every fixed $x \in M$ there is a well-defined, $\ast$-homomorphic 
 extension $\epsilon_x$ of the counit map
 $\epsilon$ to the unital $\ast$-subalgebra $\clq^\infty_x:=\{ \alpha_r(f)(x):~f \in C^\infty(M)\}$ satisfying $\epsilon_x(\alpha(f)(x))=f(x)$, where $\alpha_r$ is 
  the reduced action discussed earlier.\elmma
 {\it Proof:}\\
 Replacing $\clq$ by $\clq_r$ we can assume without loss of generality that $\clq$ has faithful Haar state and $\alpha=\alpha_r$. 
 In this case $\clq$ will have bounded antipode $\kappa$ (by Proposition \ref{tracial}). Let $\alpha_x: C^\infty(M) \raro \clq^\infty_x$ be the map defined by $\alpha_x(f):=\alpha(f)(x)$. 
 It is clearly continuous w.r.t. the Fr\'echet topology of $C^\infty(M)$ and hence the kernel say $\cli_x$
  is a closed ideal, so that the quotient which is isomorphic to $\clq^\infty_x$ is a nuclear space. Let us consider $\clq^\infty_x$ with this topology and then by nuclearity,
   the projective and injective tensor products with $\clq$ (viewed as a
separable Banach space, where separability follows from the fact that $\clq$
faithfully acts on the separable $C^*$ algebra 
   $C(M)$) coincide with $\clq^\infty_x \hat{\ot} \clq$ and the
multiplication map $m: \clq^\infty_x \hat{\ot} \clq \raro \clq$ is indeed
continuous. Now, observe that $\clq^\infty_x \hat{\ot} \clq $ is 
    isomorphic as a Fr\'echet algebra with the quotient of $C^\infty(M) \hat{\ot} \clq$ by the ideal ${\rm Ker}(\alpha_x \ot {\rm id})=\cli_x \hat{\ot} \clq$. 
         Moreover, it follows from the relation $\Delta \circ \alpha=(\alpha \ot {\rm id})\circ \alpha$ that $\Delta$ maps $\cli_x$ to $\cli_x \hat{ \ot} \clq$, 
    and in fact  it is the restriction of the Fr\'echet-continuous map $\alpha  $ there, hence induces a continuous  map from $\clq^\infty_x \cong C^\infty(M)/\cli_x$
 to $\clq^\infty_x \hat{\ot} \clq \cong (C^\infty(M) \hat{\ot} \clq)/(\cli_x \hat{\ot} \clq)$. Thus, the composite map 
 $m \circ ({\rm id} \ot \kappa) \circ \Delta: \clq^\infty_x \raro \clq$ is continuous and this coincides with 
 $\epsilon(\cdot) 1_\clq$ on the Fr\'echet-dense subalgebra of $\clq^\infty_x$ spanned by elements of the form $\alpha(f)(x)$, 
  with $f$ varying in the Fr\'echet-dense spectral subalgebra of $C^\infty(M)$. This completes the proof of the lemma.\qed\\
  \bcrlre
  \label{eps_ext}
   There is a well-defined extension of $\epsilon_x$, say $\widetilde{\epsilon_x},$  to the linear subspace spanned by elements of the 
   form $q_0q$ where $q \in \clq^\infty_x$ and $q_0 \in \clq_0$, given by $\widetilde{\epsilon_x}(q_0q)=\epsilon(q_0) \epsilon_x(q)$. A similar conclusion will hold for subspace spanned by $\{qq_{0}, q\in \clq^{\infty}_{x}, q_{0}\in\clq_{0}\}$. 
      \ecrlre
   {\it Proof:}\\ We use the notation of Lemma \ref{smooth_bdd_counit}.       
       For a finite dimensional subspace $\cld$ of $\clq_0$ denote by $\cld_1$ the subspace spanned by $({\rm id} \ot \omega)(\Delta(\cld))$
       where $\omega$ varies over the (algebraic) dual of $\clq_0$. Thus, $\Delta(\cld) \subseteq \cld_1 \ot \clq_0$.
        Let $\clw={\rm Span}(\cld_1 \clq_x^\infty)$. As $\cld_1$ is finite dimensional, $\cld_1 \ot \clq^\infty_x$ is nuclear and the multiplication map 
         from $\cld_1 \ot \clq^\infty_x$ onto $\clw$ is continuous, hence $\clw$ can be viewed as a nuclear space in the quotient 
        topology coming from $\cld_1 \ot \clq_x^\infty$.  
       As in the proof of Lemma \ref{smooth_bdd_counit} we observe that $({\rm id} \ot \kappa)\circ \Delta$ maps the subspace $\widetilde{\cld}$ spanned by 
       $\cld \clq_x^\infty$ to $ \clw \hat{\ot} \clq$, hence $\epsilon^{\cld}_x(\cdot) 1:=m_{\clw}\circ ({\rm id} \ot \kappa)\circ \Delta$ 
       defines a continuous extension of the counit  $\epsilon$ on $\widetilde{\cld}$ w.r.t. the (nuclear) quotient Fr\'echet topology of $\widetilde{\cld}$ coming from 
        $\cld \ot \clq^\infty_x$. Here, $m_\clw$ denotes the (continuous) multiplication map from $\clw \hat{\ot} \clq$ to $\clq$.
        
        Moreover, $\widetilde{\cld} \bigcap \clq_0$ is dense in $\widetilde{\cld}$ as it contains elements of the form $q_0 \alpha(f)(x)$ for 
         $q \in \cld$ and $f$ in the (dense) spectral subalgebra of $C^\infty(M)$ and clearly, $\epsilon^\cld_x$ agrees with $\epsilon$ (the original counit defined on $\clq_0$ )
          on this dense subspace. Indeed, $\epsilon(q_0\alpha(f)(x))=\epsilon(q_0)f(x)=\epsilon(q_0)\epsilon_x(\alpha(f)(x))$ for $q_0,f$ as above, and hence 
           by continuity, we can conclude that $\widetilde{\epsilon_x}$ is the unique extension of $\epsilon$ to ${\rm Span}(\cld \clq^\infty_x)$. From the uniqueness, we see that 
            $\epsilon^\cld_x$ agrees with $\epsilon^\clu_x$ on ${\rm Span}((\cld \bigcap \clu)\clq^\infty_x)$, for any two finite dimensional subspaces $\cld, \clu $ of $\clq_0$.
            In other words, the definition of $\epsilon^\cld_x$ is independent of $\cld$, and this  gives a well-defined linear map $\widetilde{\epsilon_x}$ on the whole of $\clq_0\clq^\infty_x$, which 
             satisfies $\widetilde{\epsilon_x}(q_0q)=\epsilon(q_0)\epsilon_x(q)=\epsilon(q_0)f(x)$ for all $q_0 \in \clq$ and $q=\alpha(f)(x)$ for $f$ in the spectral subalgebra, hence 
             by continuity for all $f \in C^\infty(M)$, completing the proof of the corollary.\\
             \indent For ${\rm span} \clq_{x}^{\infty}\clq_{0}$ the proof is similar and hence omitted.
           \qed\\
  \bcrlre
  \label{automatic_unitary_impl}
  For any smooth action $\alpha$ on $C^\infty(M)$, the conditions of Theorem \ref{unitary_impl} are satisfied.
  \ecrlre
  {\it Proof:}\\
  Replacing $\clq$ by the Woronowicz subalgebra generated by $\{ \alpha(f)(x), f \in C(M), x \in M\}$ we may assume that 
   $\alpha$ is faithful. If $\alpha_r(f)=0$ for $f \in C^\infty(M)$  then by Lemma \ref{smooth_bdd_counit} applying the extended $\epsilon$ we conclude 
    $f=0$. Now, consider any positive Borel measure $\mu$ of full support on $M$, with $\phi_\mu$ being the positive functional obtained by integration w.r.t $\mu$. 
    Let $\psi:=(\phi_\mu \ot h) \circ \alpha_r$ be the positive functional which is clearly $\alpha_r$-invariant and faithful on $C^\infty(M)$, 
    i.e. $\psi(f)=0, f \in C^\infty(M)$ and $f$ nonnegative  implies $f=0$. But then by Riesz Representation Theorem there is a positive Borel measure
     $\nu$ such that $\psi(f)=\int_M f d\nu$. It follows that $\nu$ has full support, hence $\psi$ is faithful also on $C(M)$. 
     Indeed, for any nonempty open subset $U$ of $M$ there is a nonzero positive $f \in C^\infty(M)$, with $0 \leq f \leq 1,$ and support of $f$ is contained in $ U$.   
     By faithfulness of $\psi$ on $C^\infty(M)$ we get $0 < \psi(f) =\int_U f d\nu \leq \nu(U).$\qed

\subsection{Defining $d\alpha$ for a smooth action $\alpha$}
Let $\alpha$ be a smooth action of a CQG $\clq$ on a manifold $M$. Recall the $C^{\infty}(M)\hat{\ot}\clq$ bimodule $\Omega^{1}(C^{\infty}(M))\bar{\ot}\clq$. However there is also a $C^{\infty}(M)$ bimodule structure of $\Omega^{1}(C^{\infty}(M))\bar{\ot}\clq$ given by
\begin{displaymath}
 f.\Omega:=\alpha(f)\Omega, \ \Omega.f:=\Omega\alpha(f),
\end{displaymath}
for $\Omega\in\Omega^{1}(C^{\infty}(M))\bar{\ot}\clq$ and $f\in C^{\infty}(M)$ where $\alpha(f)\Omega$ and $\Omega\alpha(f)$ denote the usual left and right $C^{\infty}(M)\hat{\ot}\clq$-bimodule multiplication.\\
\indent It is easy to identify elements $\Omega\in\Omega^{1}(C^{\infty}(M))\bar{\ot}\clq$ with $\clq$-valued smooth one form, i.e. $\Omega:M\raro \cup_{m\in M}(T^{\ast}_{m}M)\ot\clq$, such that for all $m\in M$, $\Omega(m)\in T^{\ast}_{m}M\ot\clq$ and for any coordinate neighborhood $U$ and the local coordinates $(x_{1},...,x_{n})$ around $m\in M$ we can find $\Omega_{i}\in C^{\infty}(M)\hat{\ot}\clq$, $i=1,...,n$, such that $\Omega(x)=\sum_{i=1}^{n}dx_{i}(x)\ot\Omega_{i}(x)$ for all $x\in U$. We shall usually write $dx_{i}(x)\ot\Omega_{i}(x)$ as $dx_{i}(x)\Omega_{i}(x)$ and $\Omega=\sum_{i=1}^{n}dx_{i}\Omega_{i}$ on $U$.\\
\indent This allows us to define $\widetilde{d}\equiv (d\ot {\rm id})$ from $C^{\infty}(M)\hat{\ot}\clq$ to $\Omega\in\Omega^{1}(C^{\infty}(M))\bar{\ot}\clq$ given by
\begin{displaymath}
 (\widetilde{d}F)(m):= \sum_{i=1}^{n} dx_{i}(m)(\frac{\partial F}{\partial x_{i}})(m),
\end{displaymath}
for $m\in M$ and for any local coordinate chart $(U,(x_{1},...,x_{n}))$ around $m$. Clearly this is uniquely defined by the condition 
\begin{displaymath}
 ({\rm id}\ot\omega)(\widetilde{d}F)=d(({\rm id}\ot\omega)(F))
\end{displaymath}
for all bounded linear functional $\omega$ on $\clq$. Thus $\widetilde{d}F$ does not depend on the choice of the local coordinates.\\
\indent We now have
\bthm
\label{def_dalpha}
The following are equivalent:\\
(i) There is a well defined, Fr\'echet continuous map $\beta:\Omega^{1}(C^{\infty}(M))\raro \Omega^{1}(C^{\infty}(M))\bar{\ot}\clq$ which is a $C^{\infty}(M)$ bimodule morphism, i.e. $\beta(f\omega)=\alpha(f)\beta(\omega)$ and $\beta(\omega f)=\beta(\omega)\alpha(f)$ for all $\omega\in \Omega^{1}(C^{\infty}(M))$ and $f\in C^{\infty}(M)$ and also $\beta(df)=(d\ot {\rm id})(\alpha(f))$ for all $f\in C^{\infty}(M)$.\\ 
(ii) For all $f,g\in C^{\infty}(M)$ and all smooth vector fields $\nu$ on $M$,
\begin{eqnarray}
\label{dalpha_eq}
(\nu\ot {\rm id})\alpha(f)\alpha(g)= \alpha(g) (\nu\ot {\rm id})\alpha(f) 
\end{eqnarray}

\ethm
{\it Proof:}\\
Proof of necessity:\\
We have $\beta (df.g)= (d\ot {\rm id})(\alpha(f)).\alpha(g)$, $\beta (g.df)=
\alpha(g).(d\ot {\rm id})(\alpha(f))$. But $df.g= g.df$ in $\Omega^{1}(C^{\infty}(M))$,
which gives
$(d\ot {\rm id})(\alpha(f)).\alpha(g)= \alpha(g).(d\ot {\rm id})(\alpha(f)), \ \forall
f,g\in C^{\infty}(M)$.\\
\indent Observe that as $\nu$ is a smooth vector field, $\nu$ is a Fr\'echet
continuous map
from $C^{\infty}(M)$ to $C^{\infty}(M)$. Thus it is enough to prove (2) for
$f,g$ belonging to the Fr\'echet dense subalgebra $\cla$ as in Theorem
\ref{suff_smooth}
. Let $\alpha(f)=f_{(0)}\ot f_{(1)}$ and $\alpha(g)=g_{(0)}\ot
g_{(1)}$(Sweedler's
notation).
Let $x\in M$ and $(U,x_{1},..,x_{n})$ be a coordinate neighborhood around
$x$. Then $[(d\ot {\rm id})(\alpha(f))\alpha(g)](x)= \sum_{i=1}^{n}g_{(0)}(x)\frac{\partial
f_{(0)}}{\partial x_{i}}(x)dx_{i}|_{x}f_{(1)}g_{(1)} $.\\
The condition $[(d\ot {\rm id})(\alpha(f))\alpha(g)](x)=[\alpha(g)(d\ot {\rm id})(\alpha(f))](x)$ gives
\begin{eqnarray}
 g_{(0)}(x)\frac{\partial f_{(0)}}{\partial x_{i}}(x)
f_{(1)}g_{(1)}= g_{(0)}(x)\frac{\partial f_{(0)}}{\partial x_{i}}(x)
g_{(1)}f_{(1)}
\end{eqnarray}
for all $i=1,...,n$. Now let $a_{i}\in C^{\infty}(M)$ for $i=1,...,n$ such that
$\nu(x)=\sum_{i=1}^{n}a_{i}(x)\frac{\partial}{\partial{x_{i}}}|_{x}$ for all
$x\in U$.\\
So 
\begin{eqnarray*}
&&[(\nu\ot {\rm id})(\alpha(f))\alpha(g)](x)\\
&=& \sum_{i=1}^{n}a_{i}(x)\frac{\partial f_{(0)}}{\partial x_{i}}(x)g_{(0)}(x)
f_{(1)}g_{(1)}\\
\end{eqnarray*}
and
\begin{eqnarray*}
&&[\alpha(g)(\nu\ot {\rm id})(\alpha(f))](x)\\
&=& \sum_{i=1}^{n}a_{i}(x)\frac{\partial f_{(0)}}{\partial x_{i}}(x)g_{(0)}(x)
g_{(1)}f_{(1)}\\
\end{eqnarray*}
Hence by (3) $[\alpha(g)(\nu\ot {\rm id})(\alpha(f))](x)=[(\nu\ot
{\rm id})(\alpha(f))\alpha(g)](x)$ for all $x\in M$\\
i.e. $[\alpha(g)(\nu\ot {\rm id})(\alpha(f))]=[(\nu\ot {\rm id})(\alpha(f))\alpha(g)]$ for all
$f,g\in \cla$.

\vspace{0.025in}
Proof of sufficiency:\\
This needs a number of intermediate lemmas. Let
$x\in M$ and $(U,x_{1},...,x_{n})$ be a coordinate neighborhood around
it. Choose smooth vector fields $\nu_{i}$'s on $M$ which are
$\frac{\partial}{\partial x_{i}}$ on $U$. So $[\alpha(g)(\nu_{i}\ot {\rm id})\alpha(f)](x)=\frac{\partial f_{(0)}}{\partial
x_{i}}(x)g_{(0)}(x)g_{(1)}f_{(1)}$ and $[(\nu_{i}\ot
{\rm id})(\alpha(f))\alpha(g)](x)=\frac{\partial f_{(0)}}{\partial
x_{i}}(x)g_{(0)}(x)f_{(1)}g_{(1)}$. Hence by the assumption 
$ \sum_{i}\frac{\partial f_{(0)}}{\partial
x_{i}}(x)g_{(0)}(x)dx_{i}|_{x}g_{(1)}f_{(1)}=\sum_{i}\frac{\partial
f_{(0)}}{\partial x_{i}}(x)g_{(0)}(x)dx_{i}|_{x}f_{(1)}g_{(1)}$, hence
$ [(d\ot {\rm id})(\alpha(f))\alpha(g)](x)=[\alpha(g)(d\ot
{\rm id})(\alpha(f))](x)$.\\
\indent Since $x$ is arbitrary, we conclude that $[\alpha(g)(d\ot {\rm id})(\alpha(f))]=[(d\ot
{\rm id})(\alpha(f))\alpha(g)]$ for all $f,g\in \cla$. So by Fr\'echet continuity of $d$
and $\alpha$ we can prove the result for $f,g \in
C^{\infty}(M)$.\qed\\
%
%
We use the commutativity to deduce the following:
\blmma
\label{new}
For $F\in C^{\infty}(\mathbb{R}^{n})$ and $g_{1},g_{2},..,g_{n} \in
C^{\infty}(M)$
\begin{eqnarray}
(d\ot {\rm id})\alpha (F(g_{1},..., g_{n}))= \sum_{i=1}^{n}\alpha
(\partial_{i}F(g_{1},..., g_{n}))(d\ot {\rm id})(\alpha(g_{i})),
\end{eqnarray}
where $\partial_{i}F$ denotes the partial derivative of $F$ with respect to the
ith coordinate of $\mathbb{R}^{n}$.
\elmma
{\it Proof}:\\
As $\{(g_{1}(x)\ldots g_{n}(x))|x\in M\}$ is a compact subset of
$\mathbb{R}^{n}$, for $F\in C^{\infty}(\mathbb{R}^{n})$, we get a sequence of
polynomials $P_{m}$ in $\mathbb{R}^{n}$ such that $P_{m}(g_{1},..., g_{n})$
converges to $F(g_{1},...,g_{n})$ in the Fr\'echet topology of
$C^{\infty}(M)$.\\
We see that for $P_{m}$, 
\begin{eqnarray*}
&& (d\ot {\rm id})(\alpha(P_{m}(g_{1},..g_{n})))\\
&=& (d\ot {\rm id})(P_{m}(\alpha(g_{1},...,g_{n})))\\
&=&  \sum_{i=1}^{n}\alpha (\partial_{i}P_{m}(g_{1},..., g_{n}))(d\ot
{\rm id})(\alpha(g_{i})), 
\end{eqnarray*}
using $ (d\ot {\rm id})(\alpha(f))\alpha(g)=\alpha(g)(d\ot {\rm id})\alpha(f)$ as well as the 
Leibniz rule  for $ (d\ot {\rm id}).$ \\
The lemma now follows from Fr\'echet continuity of $\alpha$ and $(d\ot
{\rm id})$.\qed

\blmma
\label{support}
Let $U$ be a coordinate neighborhood. Also let $g_{1}, g_{2},..., g_{n}\in
C^{\infty}(M)$ be such that $(g_{1}|_{U},\ldots g_{n}|_{U})$ gives a local
coordinate system on $U$. 
Then
\begin{eqnarray*}
(d\ot {\rm id})(\alpha(f))=\sum_{j=1}^{n}\alpha(\partial_{g_{j}}f)(d\ot
{\rm id})(\alpha(g_{j})),
\end{eqnarray*}
for all $f\in C^{\infty}(M)$ supported in $U$.
\elmma
{\it Proof}:\\
Let $F\in C^{\infty}(\mathbb{R}^{n})\raro \mathbb{R}$ be a smooth function such
that $f(m)= F(g_{1}(m),....,g_{n}(m)) \ \forall m\in U.$ Choose $\chi\in
C^{\infty}(M)$ with $\chi\equiv 1$ on $K=supp(f)$ and
$supp(\chi)\subset U$. Then $\chi f=f$ as $\chi\equiv 1$ on $K$. Hence
$\chi F(g_{1},...,g_{n})= f ( \chi F=\chi f= f \ on \ U , \ \chi
F=0 \
outside \ U)$. Also $\chi^{2}F(g_{1},...,g_{n})= \chi F(g_{1},...,g_{n})$, since
on $K$,
$\chi^{2}=\chi=1$ and outside $K$, $\chi^{2}F(g_{1},...,g_{n})= \chi
F(g_{1},...,g_{n})= 0$. Let $T:=\alpha(\chi)$ and
$S:=\alpha(F(g_{1},...,g_{n}))$. Also denote $(d\ot
{\rm id})\alpha(F(g_{1},...,g_{n}))$ by $S^{'}$ and $(d\ot {\rm id})(\alpha(\chi))$ by
$T^{'}$.\\
\indent So we have $T^{2}S=TS$ and by (2) we have $T^{\prime}T=TT^{\prime}$ and
$S^{\prime}S=SS^{\prime}$.

\begin{eqnarray}
T^{2}S^{\prime}&=& \alpha(\chi^{2})(d\ot {\rm id})(\alpha(F(g_{1},...,g_{n}))) \nonumber\\
&=& \alpha(\chi^{2})\sum_{i=1}^{n}\alpha (\partial_{i}F(g_{1},...,g_{n}))(d\ot
{\rm id})(\alpha(g_{i}))\nonumber \ (by \ (4))\\
&=& \alpha(\chi) \sum_{i=1}^{n} \alpha (\chi
\partial_{i}F(g_{1},...,g_{n}))(d\ot {\rm id})(\alpha(g_{i})) \nonumber \\
&=& \alpha(\chi) \sum_{i=1}^{n}\alpha(\partial_{g_{i}}f)(d\ot {\rm id})(\alpha(g_{i})) \
(as \ supp(\partial_{g_{i}}f)\subset K).
\end{eqnarray}
\begin{eqnarray}
TS^{\prime}&=& \alpha(\chi)(d\ot {\rm id})(\alpha(F(g_{1},...,g_{n}))) \nonumber\\
&=& \sum_{i=1}^{n} \alpha (\chi \partial_{i}F(g_{1},...,g_{n}))(d\ot
{\rm id})(\alpha(g_{i})) \nonumber \\ 
&=& \sum_{i=1}^{n} \alpha (\chi^{2} \partial_{i}F(g_{1},...,g_{n}))(d\ot
{\rm id})(\alpha(g_{i})) \nonumber \\
&=& \alpha(\chi)\sum_{i=1}^{n}\alpha(\partial_{g_{i}}f)(d\ot {\rm id})(\alpha(g_{i}))
\end{eqnarray}
Combining (5) and (6) we get
\begin{eqnarray}
T^{2}S^{\prime}= TS^{\prime}
\end{eqnarray}
Now
\begin{eqnarray}
&&T^{2}S = TS \nonumber\\
& \Rightarrow & (d\ot {\rm id})(T^{2}S)= (d\ot {\rm id})(TS)\nonumber \\
& \Rightarrow & 2TT^{\prime}S + T^{2}S^{\prime} = TS^{\prime} + T^{\prime}S (by \ Leibniz \ rule \ 
and \ T^{\prime}T=TT^{\prime})
\nonumber \\
& \Rightarrow & 2TT^{\prime}S= T^{\prime}S \nonumber \ (by (7)) \\
& \Rightarrow & 2\alpha(\chi)(d\ot {\rm id})(\alpha(\chi))\alpha(F(g_{1},...,g_{n})) =
(d\ot {\rm id})(\alpha(\chi))\alpha(F(g_{1},...,g_{n})) \nonumber \\
& \Rightarrow & 2\alpha(\chi^{2})(d\ot {\rm id})(\alpha(\chi))\alpha(F(g_{1},...,g_{n}))
= \alpha(\chi) (d\ot {\rm id})(\alpha(\chi))\alpha(F(g_{1},...,g_{n})) \nonumber \\  
& \Rightarrow & 2(d\ot {\rm id})(\alpha(\chi))\alpha(f) = (d\ot
{\rm id})(\alpha(\chi))\alpha(f)( \ using \ the \ assumption \ and \ \chi^{2}F=f)
\nonumber \\
&\Rightarrow & (d\ot {\rm id})(\alpha(\chi))\alpha(f) = 0
\end{eqnarray}
So 
\begin{eqnarray*}
(d\ot {\rm id})(\alpha(f)) &=& (d\ot {\rm id})(\alpha(\chi f))\\
&=& (d\ot {\rm id})(\alpha(\chi))\alpha(f) + \alpha(\chi)(d\ot {\rm id})(\alpha(f))\\
&=& \alpha(\chi)(d\ot {\rm id})(\alpha(f)) ( \ by \ (8))\\
&=& \alpha(\chi)(d\ot {\rm id})(\alpha(\chi F(g_{1},..., g_{n})))\\
&=& \alpha(\chi)(d\ot {\rm id})(\alpha(\chi)) \alpha(F(g_{1},...,g_{n}))+
\alpha(\chi^{2})(d\ot {\rm id})(\alpha(F(g_{1},...,g_{n})))\\
&=& (d\ot {\rm id})(\alpha(\chi))\alpha(f) + \alpha(\chi^{2})(d\ot
{\rm id})(\alpha(F(g_{1},...,g_{n}))) (Again \ by \ assumption)\\
&=& \alpha(\chi^{2})\sum_{i=1}^{n} \alpha(\partial_{i}F(g_{1},...,g_{n}))(d\ot
{\rm id})(\alpha(g_{i})) (by \ (4) \ and \ (8))\\
&=& \sum_{i=1}^{n} \alpha(\chi^{2}\partial_{i}F(g_{1},...,g_{n}))(d\ot
{\rm id})(\alpha(g_{i}))\\
&=& \sum_{i=1}^{n} \alpha(\partial_{g_{i}}f)(d\ot {\rm id})(\alpha(g_{i}))
\end{eqnarray*}
 \qed\\
 Now to complete the proof of the theorem, we want to first define a bimodule
morphism $\beta$ extending $d\alpha$ locally, i.e. we define $\beta_{U}(\omega)$
for any coordinate neighborhood $U$ and any smooth 1-form $\omega$ supported in
$U$ as follows:\\
 Choose $C^{\infty}$ functions $g_{1}\ldots g_{n}$ as before such that they give
a local coordinate system on $U$ and $\omega$ has the unique expression $\omega=
\sum_{j=1}^{n}\phi_{j}dg_{j}$. Then define $\beta_{U}(\omega):=
\sum_{j=1}^{n}\alpha(\phi_{j})(d\ot {\rm id})\alpha(g_{j})$. We claim that $\beta_{U}$ is independent of the choice of the coordinate
functions $(g_{1},\ldots,g_{n})$, i.e. if $(h_{1},\ldots,h_{n})$ is another such
set of coordinate functions on $U$ with $\omega=\sum_{j=1}^{n}\psi_{j}dh_{j}$
for some $\psi_{j}$'s in $C^{\infty}(M)$, then \\
\vspace{.15 in}
 $\sum_{j=1}^{n}\alpha(\phi_{j})(d\ot {\rm id})(\alpha(g_{j}))=
\sum_{j=1}^{n}\alpha(\psi_{j})(d\ot {\rm id})(\alpha(h_{j}))$.\\
\indent Indeed let $\chi$ be a smooth function which is
$1$ on the support of $\omega$ and $0$ outside $U$. We have
$F_{1},\ldots,F_{n}\in C^{\infty}(\mathbb{R}^{N})$ such that
$g_{j}=F_{j}(h_{1},\ldots, h_{n})$ for all $j=1,\ldots,n$ on $U$. Then $\chi
g_{j}=\chi F_{j}(h_{1},\ldots, h_{n})$ for all $j=1,\ldots,n$ . Hence $dg_{j}=\sum_{k=1}^{n}\partial_{h_{k}}(F_{j}(h_{1},\ldots,h_{n}))dh_{k}$
on $U$. That is
$\omega=\sum_{j,k}\chi\phi_{j}\partial_{h_{k}}(F_{j}(h_{1},\ldots,h_{n}))dh_{k}
$. So
$\psi_{k}=\sum_{j}\chi\phi_{j}\partial_{h_{k}}(F_{j}(h_{1},\ldots,h_{n}))$.\\
 Also, note that, as $\chi\equiv 1$ on the support of $\phi_{j}$ for all $j$, we
must have $\phi_{j}\partial_{h_{k}}(\chi)\equiv 0$, so $\chi
\phi_{j}\partial_{h_{k}}(F_{j}(h_{1},\ldots,h_{n}))= \chi\phi_{j}
\partial_{h_{k}}(\chi F_{j}(h_{1},\ldots,h_{n}))$. Thus
 
 \begin{eqnarray*}
 && \sum_{k}\alpha(\psi_{k})(d\ot {\rm id})(\alpha(h_{k}))\\
 &=&
\sum_{k,j}\alpha(\chi\phi_{j}\partial_{h_{k}}(F_{j}(h_{1},\ldots,h_{n})))(d\ot
{\rm id})(\alpha(h_{k}))\\
 &=& \sum_{k,j}\alpha(\phi_{j})\alpha (\partial_{h_{k}}(\chi
F_{j}(h_{1},\ldots,h_{n})))(d\ot {\rm id})(\alpha(h_{k}))\\
 &=& \sum_{j}\alpha(\phi_{j})(d\ot {\rm id})(\alpha(\chi F_{j}(h_{1},\ldots,h_{n}))) \
(by \ Lemma \ \ref{new})\\
 &=& \sum_{j}\alpha(\phi_{j})(d\ot {\rm id})(\alpha(\chi g_{j}))\\
 &=& \sum_{j}\alpha(\phi_{j})(d\ot {\rm id})(\alpha(g_{j}))\\
 \end{eqnarray*}
 Where the last step follows from Leibniz rule and the fact that 
 \begin{eqnarray*}
 &&\alpha(\phi_{j})(d\ot {\rm id})(\alpha (\chi))\\
 &=&\sum_{k} \alpha(\phi_{j})\alpha(\partial_{h_{k}}(\chi))(d\ot {\rm id})(\alpha
(h_{k}))\\
 &=& \sum_{k}\alpha(\phi_{j}\partial_{h_{k}}(\chi))(d\ot {\rm id})(\alpha (h_{k}))\\
 &=&0 \ (using \ \phi_{j}\partial_{h_{k}}(\chi))\equiv0),\\
 \end{eqnarray*}
 which proves the claim.\\
 \indent Hence the definition is indeed independent of choice of coordinate
system. Then for any two coordinate neighborhoods $U$ and $V$,
$\beta_{U}(\omega)=\beta_{V}(\omega)$ for any $\omega$ supported in $U\cap V$.
It also follows from the definition and the given condition (2) that $\beta_{U}$ is a $C^{\infty}(M)$-bimodule morphism. Moreover we get from Lemma \ref{support} that  $\beta_{U}(df)=(d\ot {\rm id})\alpha(f)$ for
all $f\in C^{\infty}(M)$ supported in $U$. Now we define $\beta$ globally as
follows:\\
\indent Choose (and fix) a smooth partition of unity
$\{\chi_{1},\ldots,\chi_{l}\}$
subordinate to a cover $\{U_{1},\ldots,U_{l}\}$ of the manifold $M$ such that
each $U_{i}$ is a coordinate neighborhood. Define $\beta$ by:
 \begin{displaymath}
 \beta(\omega):= \sum_{i=1}^{l}\beta_{U_{i}}(\chi_{i}\omega), 
 \end{displaymath}
 for any smooth one form $\omega$. Then for any $f\in C^{\infty}(M)$,
 \begin{eqnarray*}
 \beta(df)&=& \sum_{i=1}^{l} \beta_{U_{i}}(\chi_{i}df)\\
 &=& \sum_{i=1}^{l} \beta_{U_{i}}(d(\chi_{i}f)-fd\chi_{i})\\
 &=& \sum_{i=1}^{l}[(d\ot {\rm id})(\alpha(\chi_{i}f))-\alpha(f)(d\ot
{\rm id})(\alpha(\chi_{i}))]\\
 &=& \sum_{i=1}^{l} \alpha(\chi_{i})(d\ot {\rm id})(\alpha(f)) \ (by \ Leibniz \
rule)\\
 &=& (d\ot {\rm id})(\alpha(f))\\
 \end{eqnarray*}
 This completes the proof of the Theorem \ref{def_dalpha}.\\
  \qed\\  
  
\indent We end this subsection with an example of Hopf-algebra (of non compact type) having coaction on a coordinate algebra of an algberaic variety which violates the condition (\ref{dalpha_eq}). However as mentioned in the introduction, we don't have any example of a smooth CQG action which violates the condition (\ref{dalpha_eq}).\\
\indent Let $\cla\equiv \mathbb{C}[x]$ be the * algebra of polynomials
in one variable with complex coefficients and $\clq_{0}$ be the Hopf *
algebra generated by $a, a^{-1}, b$ subject to the following relations 
$$aa^{-1}=a^{-1}a=I, ab=q^{2}ba,$$
where $q$ is a parameter as described in \cite{loc_comp}. This Hopf algebra corresponds to the quantum $ax+b$ group. There are at least two different  (non-isomorphic)
 constructions of the analytic versions of this quantum group i.e. as locally compact quantum groups
 in the sense of \cite{loc_cpt_vaes}, one by Woronowicz (\cite{loc_comp}) the other by Baaj-Skandalis (\cite{Bj-Sk}) and Vaes-Veinermann (\cite{vaes_vein}). The coproduct is given
by (see
\cite{loc_comp} for details)
$$\Delta(a)=a\ot a, \Delta(b)=a\ot b+ b\ot I .$$
We have a coaction $\alpha:k[x]\raro k[x]\ot \clq_{0}$ given by $$\alpha(x)=x\ot
a + 1\ot b.$$
The algebra $\cla$ is the algebraic geometric analogue of
$C^{\infty}(\mathbb{R})$ and we have the following canonical derivation
$\delta:\cla\raro \cla$ corresponding to the vector field $\frac{d}{dt}$ of
$\mathbb{R}$:
$$\delta(p)=p^{\prime},$$
where $p^{\prime}$ denotes the usual derivative of the polynomial $p$.\\
However an easy computation gives
$$(\delta \ot {\rm id})\alpha(x)=1\ot a,$$
which do not commute with $\alpha(x)$ as $ab\neq ba$.

\section{Actions which preserve some Riemannian inner product}

\bdfn
\label{def_inn}
Suppose that $M$ has a Riemannian structure with the corresponding  $C^{\infty}(M)$ valued inner product $<<\cdot, \cdot>>$ on $\Omega^1(C^\infty(M))$. We call a smooth action $\alpha$ on $M$ to be inner product preserving on a Fr\'echet-dense unital $\ast$-subalgebra $\cla$ of 
$C^\infty(M)$ if 
\begin{eqnarray}
\label{inner_prod_pres_111}
 <<(d\ot {\rm id})\alpha (f),(d\ot {\rm id})\alpha (g)>>=\alpha(<<df,dg>> )
 \end{eqnarray}
 for all $f,g\in \cla$. When $\cla=C^{\infty}(M)$ we simply call the smooth action simply inner product preserving. 
\edfn
Given a smooth  action $\alpha$ of a CQG $\clq$ on a Riemannian manifold $M$,  it is easy to see, by Fr\'echet continuity of the maps $d$ and
$\alpha$ that the equation (\ref{inner_prod_pres_111}) with $f,g$ varying in any
Fr\'echet dense $\ast$-subalgebra of $C^{\infty}(M)$ is equivalent to having it for all $f,g \in C^\infty(M)$.

 
\bthm
If $\alpha$ is a smooth action of a CQG $\clq$ on a compact, smooth manifold. Then it admits a lift $d\alpha$ as a well defined $\alpha$-equivariant unitary representation on the Hilbert bimodule of one forms satisfying $d\alpha(df)=(d\ot{\rm id})\alpha(f)$ if and only if there is a Riemannian structure on the manifold such that $\alpha$ is inner product preserving.
\ethm
{\it Proof}:\\
First we suppose that $\alpha$ is inner product preserving with respect to some Riemannian inner product. Let $\cla$ be the maximal $\ast$-subalgebra of $C^{\infty}(M)$ over which the action is algebraic. Then by the remark after definition (\ref{def_inn}), we see that $\alpha$ is inner product preserving on $\cla$. For $\omega=\sum_{i}f_{i}dg_{i}$, for $f_{i},g_{i}\in \cla$ we
define
$d\alpha(\omega):=\sum_{i}(d\ot {\rm id})\alpha(g_{i})\alpha(f_i)$. Also let
$\eta=\sum_{i}f_{i^{\prime}}dg_{i^{\prime}}$,
where $f_{i^{\prime}},g_{i^{\prime}}\in \cla$,
then  $<<d\alpha(\omega),d\alpha(\eta)>> =\alpha(<<\omega,\eta>>)$. Then it
is a well defined $\alpha$ equivariant bimodule morphism. The coassociativity
condition
follows from that of $\alpha$. Moreover, as Sp $\alpha(\cla)(1\ot
\clq_{0})=\cla\ot \clq_{0}$, we have  ${\rm Sp} \ d\alpha(\Omega^{1}(\cla))(1\ot
\clq_{0})=\Omega^{1}(\cla)\ot \clq_{0}$.\\
\indent Now we prove the converse. This is actually an adaptation of the proof of Theorem 3.1 of \cite{average_1}. In this case we no longer need the assumption I of that theorem thanks to Lemma \ref{smooth_bdd_counit} (also see \cite{rigid_1}). We can assume without loss of generality that the quantum group is reduced, i.e. $\clq=\clq_r$, hence  the Haar state $h$ is faithful and the antipode is bounded.
  To justify this, observe that as $\alpha$ is algebraic on $\cla$, the
corresponding reduced action $\alpha_{r}$ (say) of  $\clq_r$ given by
$({\rm id} \ot \pi_r)\circ \alpha$ (where $\pi_r$
  denotes the quotient morphism from $\clq$ onto $\clq_r$) coincides with $\alpha$ on $\cla$. For this reason, it is enough to show that the reduced action $\alpha_r$ is inner product preserving 
   for some Riemannian structure on $\cla$. In other words, we can replace $\clq$ by $\clq_r$ in this proof, and as we already know that $\clq_r$ must be a Kac algebra, 
   antipode is norm-bounded on it.
   This proves the claim.
   
  Now, we briefly give the construction of a Riemannian metric which will be preserved by the action, following the lines of arguments in \cite{average_1}. Basically, we have to note 
   that one does not need Assumption I of that paper. Choose and fix a Riemannian inner product on $M$, with $<<\cdot, \cdot>>$ denoting the 
    corresponding $C^\infty(M)$-valued inner product on one-forms and as before, consider
     the Fr\'echet dense $\ast$-algebra $\cla$ on which $\alpha$ is algebraic. However, unlike \cite{average_1},
     $<< df, dg>>$ may  not belong to $\cla$ for $f,g \in \cla.$ However, we'll see the arguments of \cite{average_1} still work.  Extend the definition of the map $\Psi$ defined in 
    Lemma 3.3 of that paper to $C^\infty(M) \ot \clq_0$ (algebraic tensor product), i.e. let $\Psi(F)$ in $C^\infty(M)$ be given by 
    $$ \Psi(F)(x):= h (m \circ  (\kappa\ot {\rm id})(G)),$$  where $G=(\alpha\ot
{\rm id})(F)(x)$ and 
  where $m$, $\kappa$ etc. are as in that paper, $h$ the (faithful) Haar state of $\cla$. Note that $m$ is well-defined on $\clq \ot \clq_0$, so there is no problem in the above definition.
    At this point, let us make a useful observation: in the definition of $\Psi$, we may put $\kappa$ on the rightmost position, i.e. 
    we have $\Psi(F)(x)=\Psi^\sharp (F)(x):= h (m \circ  ({\rm id}\ot\kappa)(G))$. By continuity of the maps involved, it suffices to verify this for $F=f \ot q$, where $f \in \cla, q \in \clq_0$,
     and in this case, we have the following using $\kappa^2={\rm id},$ $h\circ \kappa=h$ and the traciality of $h$:
     $$\Psi(F)=f_{(0)}h(\kappa(f_{(1)})q)=f_{(0)}h(\kappa(q)f_{(1)})=f_{(0)}h(f_{(1)}\kappa(q))=\Psi^\sharp(F).$$
      The proof of complete positivity of $\Psi$ as in \cite{average_1} goes through for the extended map without any change in the arguments, i.e. the conclusion of Lemma 3.3 of that 
   paper holds. For
$\omega,\eta\in
\Omega^{1}(\cla)$ We define 
$$<<\omega,\eta>>^{\prime}:=\Psi(<<d\alpha(\omega),d\alpha(\eta)>>)$$  and observe that the proof of  Lemma 3.4, Lemma 3.5 and Lemma 3.6 of \cite{average_1} go through almost verbatim
 in our case. In fact, the only care necessary is about $x=<<d\phi_{(0)},d\psi_{(0)}>>$ (in the notation of \cite{average_1}) which no longer belongs to $\cla$, so writing Sweedler-type 
  notation for $\alpha(x)$ is not permitted. However, we can approximate $x$ by a sequence of elements $x_n$ in $\cla$ and do calculations similar to \cite{average_1} to get the 
   desired conclusions of Lemma 3.4, 3.5 and 3.6.  This shows $<\cdot, \cdot>^\prime$ is a nonnegative definite sesqui-linear form. 
      Let us give some details of the proof  of strict positive definiteness of $<\cdot, \cdot>^\prime$.   
          Let $v$ and $s_1, \ldots, s_n$ be as in the proof of Theorem 3.1 of \cite{average_1} and suppose  $<v,v>^{\prime}=0$ i.e. 
$$\sum_{i,j}\bar{c_{i}}c_{j}<<ds_{i},ds_{j}>>^{\prime}(x)=0,$$
where $v=\sum_{i}c_{i}ds_{i}(x)\in T_{x}^{\ast}(M)$. By faithfulness of $h$ and also using the observation $\Psi=\Psi^\sharp$, which allows to replace $\kappa \ot {\rm id}$ by $ {\rm id} \ot \kappa$, 
 we get  the following:
$$\sum_{i,j}\bar{c_{i}}c_{j}(({\rm id}\ot m)({\rm id}\ot {\rm id} \ot  \kappa )(\alpha\ot
{\rm id})<<d\alpha(ds_{i}),d\alpha(ds_{j})>>)(x)=0,$$ i.e. 
$$\sum_{i,j}\bar{c_{i}}c_{j}\alpha(<<ds_{i(0)},ds_{j(0)}>>)(x)\kappa(s_{i(1)}^*s_{j(1)})=0.$$
Applying the extension of 
$\epsilon $ constructed in Corollary \ref{eps_ext} to the above equation and using $\epsilon \circ \kappa=\epsilon$ on $\clq_0$, we get
\be \label{999900} \sum_{i,j}\bar{c_{i}}c_{j}<<ds_{i(0)},ds_{j(0)}>>(x)\epsilon(s_{i(1)}^*s_{j(1)})=0.\ee
But for $f,g \in \cla$, $<< df_{(0)},dg_{(0)}>>\epsilon(f_{(1)}^*g_{(1)})=<<df,dg>>,$ using $f_{(0)}\epsilon(f_{(1)})=f$ for all $f \in \cla$. Thus, 
 Using the fact that $\epsilon$ is $\ast$-homomorphism, (\ref{999900}) reduces to 
$$<\sum_{i}c_{i}ds_{i}(x),\sum_{i}c_{i}ds_{i}(x)>=0,$$
i.e $<v,v>=0.$
   \qed\\

\section{Appendix: Proof of\\ 
$C^{\infty}(M,\cla)\hat{\ot}\clb \cong C^\infty(M, \cla \hat{\ot} \clb)$}
\label{appendix}
Let $M$ be a smooth, compact $n$-dimensional manifold and $\cla$, $\clb$ be two
$C^{\ast}$ algebras. For any nice algebra $\clc$, let $Der(\clc)$ denote the set of its all closable derivations. 
From the definition of topological tensor product of two nice algebras in our sense, and as $\clb$ is a $C^{\ast}$ algebra, 
$C^{\infty}(M,\cla)\hat{\ot}\clb$ is the completion of
$C^{\infty}(M,\cla)\ot\clb$ with respect to the family of seminorms given by the
family of closable $\ast$ derivations of the form $\{\widetilde{\eta}:=\eta \ot {\rm id}:\eta\in Der(C^{\infty}(M,\cla))$, where $\eta$ is any (closable $\ast$)-derivation on 
 $C^\infty(M, \cla)$. On the other hand, we have the natural Fr\'echet topology on $C^\infty(M, \cla \hat{\ot} \clb)$ coming from the 
  derivations of the form $\delta \ot {\rm id}_{\cla \hat{\ot}\clb}$, where $\delta$ is any (closable $\ast$) derivation on $C^\infty(M)$, i.e. a smooth 
   vector field. We want to show that these two topologies coincide.
\blmma
\label{der}
Consider $C^{\infty}(M,\cla)$ as a  $C^{\infty}(M)$ bimodule using the
algebra inclusion $C^{\infty}(M)\cong C^{\infty}(M)\ot 1\subset
C^{\infty}(M,\cla)$. Let $D:C^{\infty}(M)\mapsto C^{\infty}(M,\cla)$ be a
derivation. Then given any coordinate neighborhood $(U,x_{1},...,x_{n})$, there
exist $a_{1},...,a_{n}\in C^{\infty}(M,\cla)$ such that for any $m\in U$,
\begin{displaymath}
D(f)(m)=\sum_{i=1}^{n}a_{i}(m)\frac{\partial f}{\partial x_{i}}(m).
\end{displaymath}
\elmma
{\it Proof}\\
It follows by standard arguments similar to those used in proving that
any derivation on $C^{\infty}(M)$ is a vector field. \qed\\
\bcrlre
\label{derivation}
Let $\eta$ be a closable $\ast$-derivation on $C^{\infty}(M,\cla)$. Then there exists a norm
bounded $\ast$-derivation $\eta^{\cla}:\cla\raro C^{\infty}(M,\cla)$ and finitely many smooth vector fields $\xi_{ij}$, 
elements $a_{ij} \in C^\infty(M, \cla)$, $i=1,\ldots, n$, $j=1, \ldots, p$ ($p$ positive integer)
such that 
\begin{eqnarray}
\eta ( F)(m)=\sum_{i,j}a_{ij}(m)(\xi_{ij} \ot {\rm id})(F)(m)+\eta^{\cla}(F(m))(m).
\end{eqnarray}
\ecrlre
{\it Proof}:\\
Choose any finite cover $U_1, \ldots U_p$ by coordinate neighborhoods and let $\chi_1, \ldots, \chi_p$ be the associated smooth partition of unity. 
Define $\eta^{\cla}(q):=\eta(1\ot q)$. As any closed $\ast$ derivation on a
$C^{\ast}$ algebra is norm bounded, we get from Lemma \ref{der}
and the observation that $\eta(f\ot q)=\eta(f\ot 1)(1\ot q)+ (f\ot 1)\eta(1\ot
q)$ the following expression of $\eta(F)(m)$ for $m$ in a coordinate neighborhood $U_j$, say, where $(x^{(j)}_1, \ldots x^{(j)}_n)$ are the corresponding local 
 coordinates:
 \begin{displaymath}
 (\eta F)(m)=\sum_{i=1}^{n}a_{ij}(m)\frac{\partial F}{\partial
x^{(j)}_{i}}(m)+\eta^{\cla}(F(m))(m),
\end{displaymath}
where $a_{ij}\in C^{\infty}(M,\cla)$. Now, the lemma follows by taking $\xi_{ij}=\chi_j \frac{\partial }{\partial
x^{(j)}_{i}}$.
\qed\\ 

\bcrlre
$C^\infty(M, \cla)$ is a nice algebra.
\ecrlre
{\it Proof:}\\ Denote by $\tau$ the topology on $C^\infty(M, \cla)$ coming from all closable $\ast$-derivations. As the usual Fr\'echet topology 
 on this space is given by derivations of the form $\delta \ot {\rm id}$ where $\delta$ is a smooth vector field on $M$, clearly $\tau$ is stronger than 
  the usual topology. Let us show the other direction. The expression of any $\eta\in Der(C^{\infty}(M,\cla))$ given any (26) of Corollary \ref{derivation} ensures that $\eta$ is continuous w.r.t. the usual Fr\'echet topology of $C^\infty(M, \clq)$. Thus,  a sequence $F_n$ of 
 $C^\infty(M, \cla)$ which is Cauchy in the usual topology of $C^\infty(M, \cla)$ will be Cauchy in the $\tau$-topology too. 
 It follows that $\tau$ is weaker than the usual Fr\'echet  topology of $C^\infty(M, \cla)$, hence 
    the two topologies coincide. 
    \qed\\
\blmma
\label{reverse}
Let $F\in C(M,\cla\hat{\ot}\clb)$ such that for all $\omega\in \clb^{\ast}$, where $\clb^{\ast}$ denotes the space of all bounded linear functionals on $\clb$,
$({\rm id}\ot {\rm id}\ot \omega)F\in C^{\infty}(M,\cla)$. Then $F\in
C^{\infty}(M,\cla\hat{\ot}\clb)$.
\elmma
{\it Proof}:\\
We first prove it when $M$ is an open subset $U$ of $\mathbb{R}^{n}$ with
compact closure (say) $K$. We denote the standard coordinates of $\mathbb{R}^{n}$
by $\{x_{1},...,x_{n}\}$. Let us choose a point
$x^{0}=(x^{0}_{1},...,x^{0}_{n})$ on the manifold and $h,h^{\prime}>0$ such
that $(x^{0}_{1},...,x^{0}_{i}+h,...,x^{0}_{n})$ and
$(x^{0}_{1},...,x^{0}_{i}+h^{\prime},...,x^{0}_{n})$ both belong to the open set $U$
for a fixed $i\in\{1,...,n\}.$ We
shall show that $\frac{\partial F}{\partial x_{i}}(x^{0})$ exists. That is, we
have to show that 
\begin{displaymath}
 \Omega^{F}(x^{0};h):=
\frac{F(x^{0}_{1},...,x^{0}_{i}+h,...,x^{0}_{n})-F(x^{0}_{1},...,x^{0}_{i},...,
x^{0}_{n})}{h}
\end{displaymath}
is Cauchy in $\cla\hat{\ot}\clb$ as $h\raro 0$.
For that first observe that $(({\rm id}\ot{\rm id}\ot \omega)F)(x)=({\rm
id}\ot \omega)(F(x))$ for all $x\in M$ and $\omega\in \clb^{\ast}$. Now 
\begin{eqnarray*}
 &&({\rm id}\ot \omega)(\Omega^{F}(x^{0};h)-\Omega^{F}(x^{0};h^{\prime}))\\
&=& \frac{h^{\prime}((({\rm id}\ot{\rm id}\ot
\omega)F)(x^{0}_{1},...,x^{0}_{i}+h,...,x^{0}_{n})- (({\rm id}\ot{\rm id}\ot
\omega)F)(x^{0}_{1},...,x^{0}_{i},...,x^{0}_{n}))}{hh^{\prime}}\\
&&- \frac{h((({\rm id}\ot{\rm id}\ot
\omega)F)(x^{0}_{1},...,x^{0}_{i}+h^{\prime},...,x^{0}_{n})- (({\rm id}\ot{\rm
id}\ot
\omega)F)(x^{0}_{1},...,x^{0}_{i},...,x^{0}_{n}))}{hh^{\prime}}\\
&=& \frac{h^{\prime}\int_{0}^{h}\frac{\partial}{\partial x_{i}}(({\rm id}\ot{\rm
id}\ot
\omega)F)(x^{0}_{1},...,x^{0}_{i}+u,...,x^{0}_{n})du}{hh^{\prime}}\\
&&-\frac{h\int_{
0}^{h^{\prime}}\frac{\partial}{\partial x_{i}}(({\rm id}\ot{\rm id}\ot
\omega)F)(x^{0}_{1},...,x^{0}_{i}+v,...,x^{0}_{n})dv}{hh^{\prime}}\\
&=&
\frac{\int_{0}^{h}\int_{0}^{h^{\prime}}dudv\int_{v}^{u}\frac{\partial^{2}}{
\partial x_{i}^{2}} (({\rm id}\ot{\rm id}\ot
\omega)F)(x^{0}_{1},...,x^{0}_{i}+s,...,x^{0}_{n})ds}{hh^{'}},\\
 \end{eqnarray*}
where all the integrals involved above are Banach space valued Bochner
integrals. Let ${\rm sup}_{x\in K}||\frac{\partial^{2}}{\partial
x_{i}^{2}}(({\rm id}\ot{\rm id}\ot\omega)F)(x)||=M_{\omega}$. Then using the
fact that for a regular Borel measure $\mu$ and a Banach space valued function $F$, $||\int
F d\mu||\leq \int ||F||d\mu$, we get
\begin{displaymath}
 ||({\rm id}\ot \omega)(\Omega^{F}(x^{0};h)-\Omega^{F}(x^{0};h^{\prime}))||\leq
M_{\omega}\epsilon,
\end{displaymath}
where $\epsilon={\rm min}\{h,h^{\prime}\}$.
Now consider the family $\beta^{\phi}_{x_{0};h,h^{\prime}}=(\phi\ot{\rm
id})(\Omega^{F}(x^{0};h)-\Omega^{F}(x^{0};h^{\prime}))$ for $\phi\in
\cla^{\ast}$ with $||\phi||\leq 1$. For $\omega\in \clb^{\ast}$,
\begin{displaymath}
 \omega(\beta^{\phi}_{x^{0};h,h^{\prime}})=(\phi\ot{\rm id})({\rm id}\ot
\omega)(\Omega^{F}(x^{0};h)-\Omega^{F}(x^{0};h^{\prime})).
\end{displaymath}
Hence we have 

$|\omega(\beta^{\phi}_{x^{0};h,h^{\prime}})|\leq ||({\rm id}\ot
\omega)(\Omega^{F}(x^{0};h)-\Omega^{F}(x^{0};h^{\prime}))||\leq M_{\omega}\epsilon$.
 By the uniform boundedness principle we get a constant $M>0$ such that
$||(\beta^{\phi}_{x^{0};h,h^{\prime}})||\leq M\epsilon$. But
$||(\Omega^{F}(x^{0};h)-\Omega^{F}(x^{0};h^{\prime}))||={\rm sup}_{||\phi||\leq
1}||\beta^{\phi}_{x^{0};h,h^{\prime}}||$. Therefore we get 
\begin{displaymath}
 ||(\Omega^{F}(x^{0};h)-\Omega^{F}(x^{0};h^{\prime}))||\leq M\epsilon \ for \
all \ h,h^{'}.
\end{displaymath}
 Hence $\Omega^{F}(x^{0};h)$ is Cauchy as $h$ goes to zero i.e. $\frac{\partial
F}{\partial x_{i}}(x^{0})$ exists. By similar arguments we can show the
existence of higher order partial derivatives. For a general smooth, compact
manifold $M$, going to the coordinate neighborhood and applying the above
result we can show that $F\in C^{\infty}(M,\cla\hat{\ot}\clb)$.\qed\\
Applying the above Lemma for $\cla=\mathbb{C}$, we get
\bcrlre
\label{one}
For $f\in C(M,\clb)$, if $({\rm id}\ot \phi)f\in C^{\infty}(M)$ for all
$\phi\in \clb^{\ast}$, then $f\in C^{\infty}(M,\clb)$.
\ecrlre
Now we are ready to prove the main result of this appendix.
\blmma
\label{apen_main}
We have the following isomorphism of Fr\'echet $\ast$-algebras:
$$C^{\infty}(M,\cla)\hat{\ot}\clb\cong C^{\infty}(M,\cla\hat{\ot}\clb).$$ 
\elmma
{\it Proof}:\\
  First we
show that
\begin{displaymath}
  C^{\infty}(M,\cla\hat{\ot}\clb)\subseteq C^{\infty}(M,\cla)\hat{\ot}\clb,
\end{displaymath}
and the inclusion map is Fr\'echet continuous. 
To prove the above inclusion it is enough to show that if a sequence in
$C^{\infty}(M)\ot\cla\ot\clb$ is Cauchy in the topology of the L.H.S., it is
also Cauchy in the topology of the R.H.S.. This follows from the descriptions of
derivations on the algebra $C^\infty(M, \cla)$ given in the Lemma \ref{derivation} and the fact that, $\clb$ being a $C^*$ algebra, 
 the topology on the right hand side is given by the derivations $\eta \ot {\rm id}$'s, $\eta\in Der(C^{\infty}(M))$ where $\eta$'s are 
  (closable $\ast$) derivations on $C^\infty(M, \cla)$.
  
  Moreover, observe that for
any $\omega\in\clb^{\ast}$ and $F\in C^{\infty}(M,\cla)\hat{\ot}\clb$,  $({\rm id}\ot \omega)F\in C^{\infty}(M,\cla)$.
 Hence by Lemma \ref{reverse},
we get $C^{\infty}(M,\cla)\hat{\ot}\clb\subseteq C^{\infty}(M,\cla\hat{\ot}\clb)
$ as well, i.e. the two spaces coincide as sets. By the closed graph theorem we
conclude that they are isomorphic as Fr\'echet spaces.\qed\\

 {\bf Acknowledgement:}

   The first author would like to thank Prof. Marc A. Rieffel for inviting him
 to visit  the Department of Mathematics of the University of California at Berkeley, 
 where some of the initial ideas and breakthroughs of the work came. Both the authors also thank Jyotishman Bhowmick for a thorough reading of the paper. 
   
\end{document}